\documentclass{amsart}
\usepackage{amsmath}
\usepackage{amssymb}
\usepackage{amsthm}

\DeclareMathOperator{\id}{id}

\DeclareMathOperator{\tr}{tr}
\DeclareMathOperator{\Vol}{Vol}

\DeclareMathOperator{\dvol}{dvol}

\DeclareMathOperator{\Ric}{Ric}

\DeclareMathOperator{\Rm}{Rm}
\DeclareMathOperator{\diam}{diam}

\DeclareMathOperator{\End}{End}

\DeclareMathOperator{\sq}{\#}
\DeclareMathOperator{\hash}{\sharp}

\DeclareMathOperator{\Met}{Met}

\newcommand{\olambda}{\overline{\lambda}}

\newcommand{\lp}{\langle}
\newcommand{\rp}{\rangle}
\newcommand{\lv}{\lvert}
\newcommand{\rv}{\rvert}
\newcommand{\lV}{\lVert}
\newcommand{\rV}{\rVert}

\newcommand{\bCP}{\mathbb{C}P}

\newcommand{\charconstant}{characteristic constant}

\newcommand{\mC}{\mathcal{C}}

\newcommand{\mS}{\mathcal{S}}

\newcommand{\mW}{\mathcal{W}}

\newcommand{\kM}{\mathfrak{M}}

\newcommand{\bN}{\mathbb{N}}

\newcommand{\bR}{\mathbb{R}}

\newcommand{\kso}{\mathfrak{so}}

\newcommand{\comment}[1]{}

\newtheorem{thm}{Theorem}[section]
\newtheorem{prop}[thm]{Proposition}
\newtheorem{lem}[thm]{Lemma}
\newtheorem{cor}[thm]{Corollary}

\theoremstyle{definition}
\newtheorem{defn}[thm]{Definition}

\theoremstyle{remark}
\newtheorem{remark}[thm]{Remark}

\numberwithin{equation}{section}

\begin{document}

\title[The Energy of a SMMS and Applications]{The energy of a smooth metric measure space and applications}
\author{Jeffrey S. Case}
\thanks{Partially supported by NSF-DMS Grant No.\ 1004394}
\address{Department of Mathematics \\ Princeton University \\ Princeton, NJ 08544}
\email{jscase@math.princeton.edu}
\date{}
\keywords{smooth metric measure space, quasi-Einstein metric, Perelman entropy, Yamabe constant, precompactness}
\subjclass[2000]{Primary 53C21; Secondary 53C25}
\begin{abstract}
We introduce and study the notion of the energy of a smooth metric measure space, which includes as special cases the Yamabe constant and Perelman's $\nu$-entropy.  We then investigate some properties the energy shares with these constants, in particular its relationship with the $\kappa$-noncollapsing property.  Finally, we use the energy to prove a precompactness theorem for the space of compact quasi-Einstein smooth metric measure spaces, in the spirit of similar results for Einstein metrics and gradient Ricci solitons.
\end{abstract}
\maketitle

\section{Introduction}
\label{sec:intro}

Smooth metric measure spaces as usually studied in the literature are triples $(M^n,g,e^{-\phi}\dvol_g)$ of a Riemannian manifold together with a smooth measure.  Their interest as geometric objects seems to originate in the early 1980s in work of Bakry and \'Emery~\cite{BakryEmery1985} on diffusion operators.  They also play an important role in understanding collapsing sequences of Riemannian manifolds (cf.\ \cite{CheegerColding1997}) and in Perelman's recent approach to the Ricci flow~\cite{Perelman1}.

An important feature of all of these perspectives is the Bakry-\'Emery Ricci tensor $\Ric_\phi^m$ defined on a smooth metric measure space together with a dimensional parameter $m$.  This tensor generalizes the Ricci tensor to this setting from the perspective of comparison geometry (cf.\ \cite{Wei_Wylie} and references therein).  It also gives rise to a natural notion of a quasi-Einstein metric as a Riemannian metric for which $\Ric_\phi^m=\lambda g$ for some smooth function $\phi$ and constants $\lambda,m$.  Important special cases of quasi-Einstein metrics are (conformally) Einstein metrics, gradient Ricci solitons, and the bases of Einstein warped products, including static metrics (cf. \cite{Case_Shu_Wei,HePetersenWylie2010} and references therein).

In the prequel~\cite{Case2010a}, the author introduced a slightly different perspective on smooth metric measure spaces, advocating that one should regard the dimensional parameter $m$ as a part of the data of a smooth metric measure space which interpolates between usual Riemannian geometry ($m=0$) and smooth metric measure spaces the way they are studied by Perelman ($m=\infty$).  This ``interpolating'' perspective is partially supported by the observation made in~\cite{Case2010a} that most known examples of quasi-Einstein metrics in fact exist as one-parameter families of quasi-Einstein metrics, parameterized by $m\in(1,\infty]$, and sometimes by $m\in[0,\infty]$.  As a notable example, the four-dimensional compact quasi-Einstein metrics on $\bCP^2\#\overline{\bCP^2}$ constructed by L\"u, Page and Pope~\cite{LPP} converge as $m\to\infty$ to the (compact) gradient Ricci soliton discovered by Koiso~\cite{Koiso1990} and Cao~\cite{Cao1996}.

Another feature of the perspective introduced in~\cite{Case2010a} is that it includes a natural notion of a conformal transformation of a smooth metric measure space.  One reason for introducing this was to understand the perspective of Chang, Gursky and Yang~\cite{CGY0} on studying smooth metric measure spaces, which differs from the perspective based on Bakry and \'Emery's work in that they search for conformal invariants associated to such spaces.  Particularly notable in this approach is that Perelman's $\lambda$-entropy can be realized as an ``infinite-dimensional limit'' of the Yamabe constant.  As shown in~\cite{Case2010a}, it turns out that the Bakry-\'Emery and Chang-Gursky-Yang perspectives are in some sense ``dual'' to one another, being related by a specific pointwise conformal transformation.

The aims of this article are two-fold.  First, we want to better understand what is the relationship between the Yamabe constant and Perelman's $\lambda$-entropy.  To that end, using the conformal transformations and the natural scalar curvature of a smooth metric measure space introduced in~\cite{Case2010a}, we will introduce natural conformal invariants on smooth metric measure spaces which include as special cases both the Yamabe constant and Perelman's $\lambda$-entropy.  Indeed, we will also be able to understand Perelman's $\nu$-entropy in this way, giving a new, non-probabilistic interpretation of this constant.  We will explore some basic properties of these constants, including their relationship to the Yamabe constant and to Sobolev-type inequalities.  However, we will not establish that minimizers of these constants exist in general, as this lies somewhat outside the main focus of this article.

Second, we want to understand to what extent the aforementioned families of quasi-Einstein metrics are ``typical.''  To that end, we will prove a precompactness theorem for the space of compact quasi-Einstein smooth metric measure spaces, generalizing similar results for Einstein metrics~\cite{Anderson1989,BKN1989,Gao1990} and for gradient Ricci solitons~\cite{Weber2008,Zhang2006,Zhang2009}.  An important feature of our result is that it is a precompactness theorem for smooth metric measure spaces for which the dimensional parameter is allowed to vary within $(1,\infty]$; in particular, our result can be interpreted as stating that, after taking a subsequence if necessary, noncollapsing sequences of compact quasi-Einstein metrics with $m\to\infty$ converge to shrinking gradient Ricci solitons.

There are two essential points which we need to address in order to establish our precompactness theorem.  First, we will need to establish estimates on both the metric $g$ and the potential $\phi$ of a quasi-Einstein smooth metric measure space $(M^n,g,e^{-\phi}\dvol)$ which are \emph{independent} of $m$; primarily this means understanding what happens for the existing estimates when $m\to\infty$, and occasionally revising them to yield limits.  Second, and more importantly, we need to rule out the possibility that our sequences of quasi-Einstein metrics collapse in a bad way.  The key ingredient here is the aforementioned generalization of the Yamabe constant: As we will see, uniform positive lower bounds on this constant for a quasi-Einstein smooth metric measure space will imply that the manifold is $\kappa$-noncollapsed in the sense of Perelman~\cite{Perelman1}, which will in turn allow us to prove the required $\varepsilon$-regularity lemma to control sectional curvatures (cf. \cite{TianViaclovsky2008,Weber2008}).

This article is organized as follows:

In Section~\ref{sec:summary}, we briefly summarize the perspective on smooth metric measure spaces introduced in~\cite{Case2010a}, focusing primarily on those ideas and results needed in the present article.

In Section~\ref{sec:defn}, we define the aforementioned Yamabe-type constants.  In fact, there will be two such definitions.  The first, which we call the weighted Yamabe constant, is the obvious analogue of the Yamabe constant.  The second, which we call the $m$-energy, is a closer analogue of Perelman's $\nu$-entropy.  As we will see, these two constants are equivalent when $m<\infty$, though they are distinct in the limit $m\to\infty$.

In Section~\ref{sec:existence_uniqueness}, we discuss some basic issues involving the existence and uniqueness questions for the $m$-energy.  These results are really analogues of simpler results pointed out by Perelman~\cite{Perelman1} for his constant $\nu(\tau)$, with the uniqueness being the natural generalization of Obata's result~\cite{Obata1971} for Yamabe metrics which are conformally Einstein.  We then use these results to derive some useful estimates for geometric quantities associated to a quasi-Einstein smooth metric measure space in terms of its $m$-energy.

In Section~\ref{sec:kappa_noncollapsing}, we discuss the relationship between lower bounds for the $m$-energy, Sobolev-type inequalities, and noncollapsing results for the underlying smooth metric measure spaces.  For motivational purposes, we include some results under the assumption that minimizers of the $m$-energy exist, and then we specialize to the case of a quasi-Einstein smooth metric measure space as needed to prove our precompactness theorem.

In Section~\ref{sec:precpt}, we state and prove the aforementioned precompactness theorem.  The bulk of this discussion centers on establishing the $\varepsilon$-regularity lemma, as the rest of the proof essentially follows in the same way as other results of this type (cf.\ \cite{Anderson1989,Weber2008}).  In particular, we will find it expedient to introduce a slight modification of the Riemann curvature tensor, which we call the weighted Weyl curvature, and will phrase our $\varepsilon$-regularity lemma in terms of this tensor.                 % Intro

% Acknowledgments
\section*{Acknowledgments}
This paper is based on the author's Ph.D.\ dissertation~\cite{Case_dissertation}, supervised by Professor Xianzhe Dai at the University of California, Santa Barbara, to whom the author owes many thanks.  Additionally, the ideas and their presentation within this paper have benefited greatly from conversations with Robert Bartnik, Rod Gover, Pengzi Miao, Yujen Shu, Guofang Wei and William Wylie.  I would also like to thank Paul Yang for comments and suggestions which helped improve the exposition in Section~\ref{sec:precpt}.

\section{A Quick Review of SMMS}
\label{sec:summary}

For the convenience of the reader, we begin by briefly summarizing the perspective on the study of SMMS using ideas from conformal geometry as introduced in~\cite{Case2010b}, placing a particular emphasis on those results which will be important in the present article.

\subsection{Smooth Metric Measure Spaces}
\label{sec:cwms}

\begin{defn}
A \emph{smooth metric measure space (SMMS)} is a four-tuple $(M^n,g,v^m\dvol,m)$ of a Riemannian manifold $(M^n,g)$ together with its Riemannian volume element $\dvol_g$, a positive function $v\in C^\infty(M)$, and a dimensional parameter $m\in\bR\cup\{\pm\infty\}$.

Equivalently, a SMMS is a four-tuple $(M^n,g,e^{-\phi}\dvol,m)$ of a Riemannian manifold $(M^n,g)$ together with a function $\phi\in C^\infty(M)$ and a dimensional parameter $m\in\bR\cup\{\pm\infty\}$.
\end{defn}

We shall frequently denote SMMS by the triple $(M^n,g,v^m\dvol)$, with the dimensional parameter $m$ implicitly specified in the expression for the measure.  Note that, when we do this, we are not ruling out the possibility $\lv m\rv=\infty$.

A primary concern in studying SMMS is finding the ``weighted'' analogues of familiar geometric notions from Riemannian geometry.  Arguably the most natural weighted object is the weighted divergence.

\begin{defn}
Given a SMMS $(M^n,g,v^m\dvol)$, the \emph{weighted divergence $\delta_\phi$} is the (negative of the) formal adjoint of the exterior derivative $d$ with respect to the measure $v^m\dvol$; i.e.\ for all $\alpha\in\Omega^k(M)$ and $\beta\in\Omega^{k+1}(M)$, at least one of which is compactly supported,
\[ \int_M \lp \alpha, d\beta \rp\, v^m\dvol = -\int_M \lp \delta_\phi\alpha,\beta\rp\,v^m\dvol . \]

The \emph{weighted Laplacian $\Delta_\phi$} is the operator $\Delta_\phi=\delta_\phi d$ defined on $C^\infty(M)$.
\end{defn}

The weighted Laplacian can be computed in terms of the usual Laplacian by $\Delta_\phi=\Delta - \lp\nabla\phi,\nabla\cdot\rp$.  A similar formula holds for the weighted divergence, and we can of course consider adjoints of more general operators between vector bundles, as will appear in Section~\ref{sec:precpt}.

In order to accomplish our goal of introducing weighted analogues of the Yamabe constant, we must first define what we mean by a conformal transformation.

\begin{defn}
Two SMMS $(M^n,g,v^m\dvol_g)$ and $(M^n,\hat g,\hat v^m\dvol_{\hat g})$ are said to be \emph{(pointwise) conformally equivalent} if there is a positive function $u\in C^\infty(M)$ such that
\begin{equation}
\label{eqn:scms}
\left(M^n,\hat g,\hat v^m\dvol_{\hat g}\right) = \left( M^n,u^{-2}g, u^{-m-n}v^m\dvol_g\right) .
\end{equation}
\end{defn}

Note that this makes sense in the cases $\lv m\rv=\infty$ by defining $u=e^{\frac{f}{m+n-2}}$ and then taking the limit.  In this case, a ``conformal transformation'' is simply a change of measure.

It is clear that~\eqref{eqn:scms} defines an equivalence relation in the space of SMMS, allowing one to talk about conformal invariants.

\begin{defn}
Let $(M^n,g,v^m\dvol_g)$ be a SMMS.  An object $T=T[g,v^m\dvol_g]$ defined on $M$ is a \emph{SMMS invariant} if for all diffeomorphisms $\psi\colon M\to M$,
\[ \psi^\ast\left(T[g,v^m\dvol_g]\right) = T[\psi^\ast g,\psi^\ast(v^m\dvol_g)] . \]
\end{defn}

\begin{defn}
Let $(M^n,g,v^m)$ be a SMMS.  A SMMS invariant $T[g,v^m\dvol_g]$ on $M$ is said to be \emph{conformally invariant of weight $w$} if for all positive $u\in C^\infty(M)$,
\[ T\left[u^{-2}g,u^{-m-n}v^m\dvol_g\right] = u^{-(m+n)w}T\left[g,v^m\dvol_g\right] . \]
\end{defn}
\subsection{Quasi-Einstein Metrics}
\label{sec:qe}

In order to discuss quasi-Einstein metrics, we first need the weighted analogues of the Ricci and scalar curvatures.

\begin{defn}
Let $(M^n,g,v^m\dvol)$ be a SMMS.  The \emph{Bakry-\'Emery Ricci tensor $\Ric_\phi^m$} is the symmetric $(0,2)$-tensor
\[ \Ric_\phi^m := \Ric - mv^{-1}\nabla^2 v = \Ric + \nabla^2\phi - \frac{1}{m}d\phi\otimes d\phi . \]
The \emph{weighted scalar curvature $R_\phi^m$} is the function
\[ R_\phi^m := R - 2mv^{-1}\Delta v - m(m-1)\lv\nabla v\rv^2 = R + 2\Delta\phi - \frac{m+1}{m}\lv\nabla\phi\rv^2 . \]
\end{defn}

\begin{defn}
A SMMS $(M^n,g,v^m\dvol_g)$ is said to be \emph{quasi-Einstein} if there is a constant $\lambda\in\bR$ such that
\begin{equation}
\label{eqn:qe_ric}
\Ric_\phi^m = \lambda g .
\end{equation}
When this is the case, we call $\lambda$ the \emph{quasi-Einstein constant}.
\end{defn}

An important result of Kim and Kim~\cite{Kim_Kim} is the following consequence of the Bianchi identity.

\begin{lem}
\label{lem:bianchi}
Let $(M^n,g,v^m\dvol)$ be a quasi-Einstein SMMS with quasi-Einstein constant $\lambda$.  Then there is a constant $\mu\in\bR$ such that
\begin{equation}
\label{eqn:bianchi}
R_\phi^m + m\mu v^{-2} = (m+n)\lambda
\end{equation}
\end{lem}

This makes sense in the limit $\lv m\rv=\infty$, where it states that there is a constant $\mu^\prime\in\bR$ such that
\[ R_\phi^\infty + 2\lambda(\phi-n) = -\mu^\prime, \]
where one requires that $\mu=\lambda+\frac{1}{m}(\mu^\prime-n\lambda)+o\left(\frac{1}{m}\right)$; see~\cite{Case2010a} for a more precise statement.

For our perspective on SMMS, it will turn out that the constant $\mu$ defined by~\eqref{eqn:bianchi} is typically more important than the constant $\lambda$.  For this reason, it is useful to give it a name.

\begin{defn}
Let $(M^n,g,v^m\dvol)$ be a quasi-Einstein SMMS with quasi-Einstein constant $\lambda$.  The \emph{\charconstant} $\mu$ is the constant such that~\eqref{eqn:bianchi} holds.
\end{defn}

Note that, by the aforementioned limiting behavior of~\eqref{eqn:bianchi}, when $\lv m\rv=\infty$ the \charconstant\ of a quasi-Einstein SMMS is exactly the quasi-Einstein constant.

To understand the weighted Yamabe problem, we need to know how the Bakry-\'Emery Ricci tensor and the weighted Yamabe constant transform under a conformal change of a SMMS.

\begin{prop}
\label{prop:curvature_conformal}
Let $(M^n,g,v^m\dvol)$ be a SMMS and let $u\in C^\infty(M)$ be a positive function.  The Bakry-\'Emery Ricci curvature $\Ric_{f,\phi}^m$ and the weighted scalar curvature $u^2R_{f,\phi}^m$ of the SMMS $(M^n,\hat g,\hat v^m\dvol_{\hat g})$ determined by $u$ via~\eqref{eqn:scms} are given by
\begin{align*}
\Ric_{f,\phi}^m & = \Ric_\phi^m + (m+n-2)u^{-1}\nabla^2 u + \left(u^{-1}\Delta_\phi u - (m+n-1)u^{-2}\lv\nabla u\rv^2\right)g \\
R_{f,\phi}^m & = R_\phi^m + 2(m+n-1)u^{-1}\Delta_\phi u - (m+n)(m+n-1)u^{-2}\lv\nabla u\rv^2 .
\end{align*}
\end{prop}

In particular, we can easily understand what it means for a SMMS to be conformally quasi-Einstein.

\begin{prop}
\label{prop:curvature_conformal_qe}
Let $(M^n,g,v^m\dvol)$ be a quasi-Einstein SMMS and suppose that $u\in C^\infty(M)$ is such that the SMMS $(M^n,\hat g,\hat v^m\dvol_{\hat g})$ determined by~\eqref{eqn:scms} is quasi-Einstein with quasi-Einstein constant $\lambda$ and \charconstant\ $\mu$.  Then
\begin{subequations}
\label{eqn:qe_uv}
\begin{align}
\label{eqn:qe_tf_uv} 0 & = \left(uv\Ric + (m+n-2)v\nabla^2u - mu\nabla^2v\right)_0 \\
\label{eqn:qe_lambda_uv} n\lambda v^2 & = (uv)^2R + (m+2n-2)uv^2\Delta u - mu^2v\Delta v \\
\notag & \quad - (m+n-1)nv^2\lv\nabla u\rv^2 + mnuv\lp\nabla u,\nabla v\rp \\
\label{eqn:qe_mu_uv} n\mu u^2 & = (uv)^2R + (m+n-2)uv^2\Delta u - (m-n)u^2v\Delta v \\
\notag & \quad - (m+n-2)nuv\lp\nabla u,\nabla v\rp + (m-1)nu^2\lv\nabla v\rv^2 ,
\end{align}
\end{subequations}
where $T_0=T-\frac{1}{n}\tr_g T\,g$ denotes the tracefree part of a section $T\in S^2T^\ast M$ and all inner products and traces on the right hand side are computed with respect to $g$.
\end{prop}

From~\eqref{eqn:qe_uv}, one sees that $\lambda$ and $\mu$ can be regarded as the squared-lengths of $u$ and $v$, in the sense that if we set $\tilde u=cu$ and $\tilde v=kv$, then~\eqref{eqn:qe_uv} holds for $\tilde u$ and $\tilde v$ by setting $\tilde\lambda=c^2\lambda$ and $\tilde\mu=k^2\mu$.  In particular, it is natural to assign to a SMMS the \charconstant\ $\mu$ as a normalization of the measure $v^m\dvol$.  One way we will do this is in the search for quasi-Einstein scales.

\begin{defn}
Let $(M^n,g,v^m\dvol)$ be a SMMS with \charconstant\ $\mu$.  A \emph{quasi-Einstein scale} is a function $u\in C^\infty(M)$ such that~\eqref{eqn:qe_uv} holds for some constant $\lambda$.  In this case, we will still call $\lambda$ the \emph{quasi-Einstein constant}.
\end{defn}

It is clear that~\eqref{eqn:qe_uv} is invariant under the transformation
\[ (u,v,\lambda,\mu,m) \mapsto (v,u,\mu,\lambda,2-m-n) . \]
This can be recast as a ``duality'' for SMMS.

\begin{cor}
\label{cor:duality}
Let $(M^n,g,v^m\dvol)$ be a SMMS with \charconstant\ $\mu$ and suppose that $u\in C^\infty(M)$ is a positive quasi-Einstein scale with quasi-Einstein constant $\lambda$.  Then the SMMS $(M^n,g,u^{2-m-n}\dvol)$ with \charconstant\ $\lambda$ is such that $v\in C^\infty(M)$ is a quasi-Einstein scale with quasi-Einstein constant $\lambda$.
\end{cor}

We can recast Corollary~\ref{cor:duality} as stating that the following are equivalent:
\begin{enumerate}
\item $(M^n,g,1^m\dvol)$ is a SMMS with \charconstant\ $\mu$ such that the function $u\in C^\infty(M)$ is a quasi-Einstein scale with quasi-Einstein constant $\lambda$.
\item $(M^n,g,u^{2-m-n}\dvol)$ is a quasi-Einstein SMMS with quasi-Einstein constant $\mu$ and \charconstant\ $\lambda$.
\end{enumerate}
It is in this sense that we will use Corollary~\ref{cor:duality} in this article; for example, it implies that if $(M^n,g,1^m\dvol)$ is a SMMS with \charconstant\ $\mu$ which admits a quasi-Einstein scale $u=e^{\frac{f}{m+n-2}}$, then
\[ \Ric + \nabla^2 f + \frac{1}{m+n-2}df\otimes df = \mu g . \]
\subsection{Some Key Results From~\cite{Case2010a}}
\label{sec:results}

One of the main results in~\cite{Case2010a} is the variational characterization of quasi-Einstein metrics with \charconstant\ $\mu$, which is accomplished using the $(m,\mu)$-energy functional.

\begin{defn}
Let $M^n$ be a smooth manifold and fix constants $m\in\bR\cup\{\pm\infty\}$ and $\mu\in\bR$.  The \emph{$(m,\mu)$-energy functional $\mW_\mu^m\colon\Met(M)\times\kM\to\bR$} is given by
\begin{subequations}
\label{eqn:mW}
\begin{equation}
\label{eqn:mW_finite}
\mW_\mu^m\left(g,v^m\dvol_g\right) = \int_M \left( R_\phi^m + m\mu v^{-2}\right) v^m\dvol_g
\end{equation}
when $\lv m\rv<\infty$ and by
\begin{equation}
\label{eqn:mW_infinite}
\mW_\mu^m\left(g,e^{-\phi}\dvol_g\right) = \int_M \left( R_\phi^m + 2\mu(\phi-n)\right) e^{-\phi}\dvol_g
\end{equation}
\end{subequations}
when $\lv m\rv=\infty$.

When the context is clear, we shall often refer to the $(m,\mu)$-energy functional as the \emph{energy functional}.
\end{defn}

In the above definition, we are denoting by $\Met(M)$ and $\kM$ the space of Riemannian metrics and smooth measures on $M$.

While the definitions given above are not continuous in the limit $\lv m\rv\to\infty$, one can introduce a natural renormalization to make this the case.

\begin{prop}
\label{prop:energy_functional_smooth}
Let $(M^n,g)$ be a Riemannian manifold and fix $\phi\in C^\infty(M)$, $\mu\in\bR$.  Then
\[ \lim_{m\to\infty} \mW_\mu^m\left(g,e^{-\phi}\dvol_g,m\right) - (m+2n)\Vol_\phi(M) = \mW_\mu^\infty\left(g,e^{-\phi}\dvol_g,\infty\right) . \]
\end{prop}

The variational characterization of quasi-Einstein metrics follows immediately from the computation for the first variation of the energy functional.

\begin{prop}
\label{prop:gen_var}
Let $M^n$ be a smooth manifold and fix constants $m\in\bR\cup\{\pm\infty\}$ and $\mu\in\bR$.  Let $(g,e^{-\phi}\dvol_g)\in\Met(M)\times\kM$ and let $(h,\psi)=(\delta g,\delta\phi)$ be a compactly supported variation of $g$ and $h$.  Then
\begin{align*}
\delta\mW_\mu^m & = -\int_M \bigg[ \lp\Ric_\phi^m - \frac{1}{2}(R_\phi^m+m\mu v^{-2})g,h\rp \\
& \qquad + \left(R_\phi^m - \frac{2}{m}\Delta_\phi\phi + (m-2)\mu v^{-2}\right)\psi \bigg] e^{-\phi}\dvol_g ,
\end{align*}
where all derivatives and inner products are taken with respect to $g$ and we define
\[ \delta\mW_\mu^m = \frac{d}{ds}\mW_\mu^m\left(g+sh,e^{-\phi+s\psi}\dvol_{g+sh}\right)\big|_{s=0} . \]
\end{prop}

\begin{prop}
\label{prop:qe_var}
Let $M^n$ be a compact manifold, fix constants $m\in\bR\cup\{\pm\infty\}$ and $\mu\in\bR$, and denote by
\[ \mC_1(m) = \left\{ (g,e^{-\phi}\dvol_g)\in\Met(M)\times\kM \colon \int_M e^{-\phi}\dvol_g = 1 \right\} . \]
Let $(g,e^{-\phi}\dvol_g)\in\mC_1(M,m)$ and consider a variation $(h,\psi)=(\delta g,\delta\phi)$ which remains in $\mC_1(m)$.  Then
\begin{align}
\label{eqn:constraint_implication} 0 & = \int \left(\psi-\frac{1}{2}\tr_g h\right)e^{-\phi}\dvol_g \\
\label{eqn:wmumvar} \delta\mW_\mu^m & = -\int\bigg[\lp\Ric_\phi^m - \frac{1}{m}\Delta_\phi\phi\,g - \mu v^{-2}g,h\rp \\
\notag & \qquad + \left(R_\phi^m-\frac{2}{m}\Delta_\phi\phi + (m-2)\mu v^{-2}\right)\left(\psi-\frac{1}{2}\tr h\right)\bigg] e^{-\phi}\dvol_g .
\end{align}

In particular, if $(g,v^m\dvol_g)\in\mC_1(M,m)$ is a critical point of the energy functional, then $(M^n,g,v^m\dvol_g)$ is a quasi-Einstein SMMS with \charconstant\ $\mu$.
\end{prop}

As we will be primarily interested in studying compact quasi-Einstein SMMS, the following result of Kim and Kim~\cite{Kim_Kim} on the signs of the quasi-Einstein constant and the \charconstant\ will be useful.

\begin{prop}
\label{prop:compact_positive_constants}
A compact quasi-Einstein SMMS is \emph{trivial} if the quasi-Einstein constant or the \charconstant\ is nonpositive.
\end{prop}

Here, we say that $(M^n,g,e^{-\phi}\dvol)$ is trivial if $\phi$ is constant.

Another useful fact about compact quasi-Einstein SMMS is the following estimate for the gradient of the quasi-Einstein scale and for the scalar curvature.

\begin{cor}
\label{cor:dil}
Let $(M^n,g,u^{2-m-n}\dvol_g)$ be a compact quasi-Einstein SMMS with $m>1$, quasi-Einstein constant $\mu$, and \charconstant\ $\lambda$.  Then
\begin{equation}
\label{eqn:dil_estimate}
\lv\nabla u\rv^2 + \frac{\lambda}{m+n-1} < \frac{\mu}{m-1}u^2
\end{equation}
and
\begin{equation}
\label{eqn:scal_bd}
-\frac{n(n-1)}{m-1}\mu < R \leq (m+2n-2)\mu .
\end{equation}
\end{cor}

These results are actually sharp, as can be seen (via Corollary~\ref{cor:duality}) by direct computation involving the positive elliptic Gaussians.

\begin{defn}
\label{defn:peg}
Fix $m\geq0$ and define $k=\sqrt{m+n-1}$.  The \emph{positive elliptic $m$-Gaussian} is the SMMS
\[ \left(S_+^n, g=dr^2\oplus\left(k\sin(\frac{r}{k})\right)^2d\theta^2, \cos^m(\frac{r}{k})\,\dvol_g\right) ,\]
where $d\theta^2$ denotes the standard metric of constant sectional curvature, so that the metric $g$ is the standard metric on $S^n$, normalized so that $\Ric_g=\frac{n-1}{m+n-1}$, and $S_+^n$ is the hemisphere $\{r<\frac{k}{2}\pi\}$.

The \emph{positive elliptic $\infty$-Gaussian} is the SMMS
\[ \left( [0,\infty)\times S^{n-1}, dr^2\oplus r^2d\theta^2, e^{-\frac{r^2}{2}}\dvol_g, \infty\right) . \]
\end{defn}

The normalizations chosen here is so that the positive elliptic Gaussians form a smooth family of quasi-Einstein SMMS for $m\in[0,\infty]$.

\begin{prop}
\label{prop:peg_properties}
The positive elliptic $m$-Gaussian is a quasi-Einstein SMMS with quasi-Einstein constant $\lambda=1$ and \charconstant\ $\mu=\frac{m-1}{m+n-1}$.  Moreover, the positive elliptic $m$-Gaussians converge to the positive elliptic $\infty$-Gaussian as $m\to\infty$ in the pointed measured Cheeger-Gromov sense, where we have fixed the base point $r=0$.
\end{prop}

As defined, the positive elliptic Gaussians are somewhat awkward to work with because they are not complete manifolds, or equivalently, because the potential $v$ degenerates along their boundary.  For this reason, we will find it useful to instead consider the conformally equivalent SMMS
\begin{equation}
\label{eqn:hyperbolic_peg}
\left( H^n, dr^2\oplus\left(k\sinh(\frac{r}{k})\right)^2d\theta^2, 1^m\dvol_g \right) ,
\end{equation}
which arises by using $v$ as the conformal factor in~\eqref{eqn:scms}; in particular, $u=\cosh(\frac{r}{k})$ is a quasi-Einstein scale for~\eqref{eqn:hyperbolic_peg}.  The benefit of considering the positive elliptic Gaussian this way is that hyperbolic space is complete, which in particular will make it easier to formulate the estimates we need to establish the precompactness theorem of Section~\ref{sec:precpt}.                % Summary of SMMS
\section{The Energy of a SMMS}
\label{sec:defn}

We have seen that quasi-Einstein metrics can be characterized as critical points of the energy functional; indeed, it is straightforward to check that, like Einstein metrics, they are saddle points.  As such, one might hope to find quasi-Einstein metrics via a minimax construction.  While this is a difficult problem in general, being that it is still open in the Einstein case, it is natural to expect that we can get some understanding of the ``mini'' part; i.e.\ we can understand the problem of minimizing the energy functional within the conformal class of a SMMS.

The purpose of this section is to introduce conformal invariants associated to a SMMS which allow us to make sense of what this means.  It will in fact be convenient to introduce two different invariants, which can be thought of as generalizations of Perelman's $\lambda$ and $\nu$-entropies~\cite{Perelman1}, though they will eventually be seen to be equivalent in a certain sense.  Once these definitions are made, we will consider the relationship between our two conformal invariants, as well as the relationship between them and the Yamabe constant and Perelman's entropies, and also establish some of their basic variational properties.  In the ensuing sections, we will investigate some of the basic properties of these invariants, especially as they will be needed to establish the precompactness theorem in Section~\ref{sec:precpt}.

For convenience, unless otherwise specified, we assume for the remainder of this article that the dimensional parameter $m$ of a SMMS is nonnegative.

\subsection{The weighted Yamabe constant}
\label{sec:yamabe}

At first glance, the difference between our two conformal invariants lies in the choice of \charconstant.  Considering the \charconstant\ to be zero yields a definition which appears as an obvious analogue of the Yamabe constant, giving it its name.

\begin{defn}
Let $(M^n,g,v^m\dvol)$ be a compact SMMS.  The \emph{weighted Yamabe constant $\sigma_{1,2}(g,v^m\dvol)$} is defined by
\begin{equation}
\label{eqn:weighted_yamabe_conformal}
\sigma_{1,2}(g,v^m\dvol) = \inf\left\{ \int_M R_\phi^m\left[\hat g,\hat v^m\dvol_{\hat g}\right]\,\hat v^m\dvol_{\hat g} \right\} ,
\end{equation}
where, given a positive function $u\in C^\infty(M)$, we define
\[ (\hat g,\hat v^m\dvol_{\hat g}) = \left(u^{-2}g,u^{-m-n}v^m\dvol_g\right), \]
and the infimum is taken over all such $(\hat g,\hat v^m\dvol)\in\mC_{1,2}(M)$, defined by
\[ \mC_{1,2}(M) = \left\{ (g,v^m\dvol_g)\in\Met(M)\times\kM \colon \int_M v^m\dvol_g = 1 = \int_M v^{m-2}\dvol_g \right\} . \]
\end{defn}

Here, one should regard the space $\mC_{1,2}(M)$ as the space of all pairs $(\hat g,\hat v^m\dvol_{\hat g})$ which are conformally equivalent to $(g,v^m\dvol_g)$ subject to a normalization condition on both $g$ and $v$.  The second integral in the definition of $\mC_{1,2}(M)$ is used, as opposed to some other integral involving different powers of $v$, because it is the natural integral appearing when one integrates the equation~\eqref{eqn:bianchi}.  It is also for this reason that we denote the weighted Yamabe constant with the subscripts $\sigma_{1,2}$: The first subscript indicates that we are considering the weighted scalar curvature, and the second subscript indicates that our normalization includes the $L^2(M,v^m\dvol_g)$-norm of $v^{-1}$.  Finally, it is clear that the weighted Yamabe constant is a conformal invariant of the SMMS $(M^n,g,v^m\dvol)$.

There are two equivalent ways to formulate the weighted Yamabe constant which are useful for different purposes.  First, we observe that since each of the integrals involved in the definition of $\sigma_{1,2}$ is homogeneous in the rescalings $g\mapsto c^2g$ and $v\mapsto kv$ for $c,k>0$, the weighted Yamabe constant is actually scale invariant.  In particular, we can reformulate the weighted Yamabe constant as a quotient.

\begin{prop}
Let $(M^n,g,v^m\dvol)$ be a compact SMMS.  The weighted Yamabe constant is equivalently defined by
\begin{equation}
\label{eqn:weighted_yamabe_conformal_homogeneous}
\sigma_{1,2}(g,v^m\dvol) = \inf\left\{ \frac{(\int R_\phi^m[\hat g,\hat v^m\dvol_{\hat g}]\hat v^m\dvol_{\hat g})\,(\int \hat v^{m-2}\dvol_{\hat g})^{m/n}}{(\int_M \hat v^m\dvol_{\hat g})^{(m+n-2)/n}} \right\} ,
\end{equation}
where again we write $(\hat g,\hat v\dvol_{\hat g})=(u^{-2}g,u^{-m-n}v^m\dvol_g)$ and now take the infimum over all positive $u\in C^\infty(M)$.
\end{prop}

\begin{proof}

It suffices to check that the quotient inside the infimum of~\eqref{eqn:weighted_yamabe_conformal_homogeneous} is invariant under the simultaneous rescalings $g\mapsto c^2g$ and $v\mapsto kv$.  To that end, we check that
\begin{align*}
\int R_\phi^m\left[c^2g,(kv)^m\dvol_{c^2g}\right] (kv)^m\dvol_{c^2g} & = c^{n-2}k^m\int R_\phi^m\left[g,v^m\dvol_g\right] v^m\dvol_g \\
\left(\int (kv)^{m-2}\dvol_{c^2g}\right)^{\frac{m}{n}} & = c^mk^{\frac{m(m-2)}{n}}\left(\int v^{m-2}\dvol_g\right)^{\frac{m}{n}} \\
\left(\int (kv)^m\dvol_{c^2g}\right)^{\frac{m+n-2}{n}} & = c^{m+n-2}k^{\frac{m(m+n-2)}{n}}\left(\int v^m\dvol_g\right)^{\frac{m+n-2}{n}} .
\end{align*}
The result follows immediately.
\end{proof}

The second way to rephrase the weighted Yamabe constant is by explicitly writing the formulae in~\eqref{eqn:weighted_yamabe_conformal} in terms of $u$.  If one does this, one finds that the weighted Yamabe constant is a natural geometric invariant associated to the weighted conformal Laplacian.

\begin{defn}
Let $(M^n,g,v^m\dvol)$ be a SMMS.  The \emph{weighted conformal Laplacian $L_\phi^m\colon C^\infty(M)\to C^\infty(M)$} is the operator
\[ L_\phi^m = -\frac{4(m+n-1)}{m+n-2}\Delta_\phi + R_\phi^m , \]
where we regard the second summand as a multiplication operator.
\end{defn}

Besides being a natural generalization of the conformal Laplacian, which is the Riemannian case $m=0$, our choice of terminology is justified by the following lemma.

\begin{lem}
\label{lem:weighted_conformal_laplacian}
Let $(M^n,g,v^m\dvol)$ be a SMMS and write $L_\phi^m[g,v^m\dvol]$ for the weighted conformal Laplacian.  Given any positive $u\in C^\infty(M)$, it holds that
\begin{equation}
\label{eqn:weighted_conformal_laplacian_covariant}
L_\phi^m\left[u^{-2}g,u^{-m-n}v^m\dvol_g\right] = u^{\frac{m+n+2}{2}} \circ L_\phi^m\left[g,v^m\dvol_g\right] \circ u^{-\frac{m+n-2}{2}},
\end{equation}
where we regard the factors $u^{\frac{m+n+2}{2}}$ and $u^{-\frac{m+n-2}{2}}$ as multiplication operators.
\end{lem}

\begin{proof}

This follows using the formula for how the Laplacian changes (cf.\ \cite{Besse}) under a conformal transformation and from Proposition~\ref{prop:curvature_conformal}.
\end{proof}

\begin{cor}
\label{cor:weighted_conformal_laplacian_scalar}
Let $(M^n,g,v^m\dvol)$ be a SMMS and let $u\in C^\infty(M)$ be a positive function.  Then
\[ R_\phi^m\left[u^{-2}g,u^{-m-n}v^m\dvol_g\right] = u^{\frac{m+n+2}{2}} L_\phi^m\left[g,v^m\dvol_g\right]\left(u^{-\frac{m+n+2}{2}}\right) . \]
\end{cor}

\begin{proof}

This follows immediately from Lemma~\ref{lem:weighted_conformal_laplacian} and the observation that
\[ R_\phi^m\left[g,v^m\dvol_g\right] = L_\phi^m\left[g,v^m\dvol_g\right](1) . \qedhere \]
\end{proof}

\begin{cor}
\label{cor:weighted_yamabe_conformal_functional}
Let $(M^n,g,v^m\dvol)$ be a SMMS.  Then the weighted Yamabe constant is equivalently defined by
\begin{equation}
\label{eqn:weighted_yamabe_conformal_functional}
\sigma_{1,2}(g,v^m\dvol) = \inf\left\{ \frac{(L_\phi^mw,w)\lV wv^{-1}\rV_2^{2m/n}}{\lV w\rV_p^q} \colon 0\not=w\in L_1^2(M,v^m\dvol) \right\},
\end{equation}
where $p=\frac{2(m+n)}{m+n-2}$, $q=\frac{2(m+n)}{n}$, $(f,g)$ denotes the usual inner product on $L^2(M,v^m\dvol)$, and all norms are computed with respect to the measure $v^m\dvol$.
\end{cor}

\begin{proof}

Given a positive function $u\in C^\infty(M)$, define $w=u^{-\frac{m+n-2}{2}}$.  It follows immediately from Corollary~\ref{cor:weighted_conformal_laplacian_scalar} that the weighted Yamabe constant is given by~\eqref{eqn:weighted_yamabe_conformal_functional} with the infimum taken over all positive $w\in C^\infty(M)$.  Since $\big\|\nabla\lv w\rv\big\|_2\leq\lV\nabla w\rV_2$ for all $w\in C^\infty(M)$, the weighted Yamabe constant is equivalently defined by taking the infimum over all nonzero $w\in C^\infty(M)$.  Finally, since $C^\infty(M)$ is dense in $L_1^2(M,v^m\dvol)$, the result follows.
\end{proof}

We now consider some basic properties of the weighted Yamabe constant.  To that end, the following straightforward consequence of H\"older's inequality and the Sobolev embedding $L_1^2(M,v^m\dvol)\hookrightarrow L^{\frac{2(m+n)}{m+n-2}}(M,v^m\dvol)$ will be useful.

\begin{lem}
\label{lem:holder}
Let $(M^n,g,v^m\dvol)$ be a compact SMMS.  For any $w\in L_1^2(M,v^m\dvol)$, it holds that
\begin{align*}
\lV w\rV_2^2 & \leq \lV w\rV_p^2 \Vol_\phi(M)^{\frac{2}{m+n}} \\
\lV wv^{-1}\rV_2^2 & \leq \lV w\rV_p^2 \left(\int v^{-n}\dvol\right)^{\frac{2}{m+n}},
\end{align*}
where $p=\frac{2(m+n)}{m+n-2}$.
\end{lem}

In particular, this allows us to show that the weighted Yamabe constant of a compact SMMS is always finite.

\begin{prop}
\label{prop:m-yamabe_finite}
Let $(M^n,g,v^m\dvol)$ be a compact SMMS.  Then there is a constant $C$ such that $\sigma_{1,2}(g,v^m\dvol)\geq C$.
\end{prop}

\begin{proof}

Let $w\in L_1^2(M,v^m\dvol)$ and set $C_1=\sup\{-R_\phi^m,0\}\geq0$.  It is clear that
\[ (L_\phi^mw,w) \geq -C_1\lV w\rV_2^2 . \]
It then follows from Lemma~\ref{lem:holder} that
\[ \sigma_{1,2}(g,v^m\dvol) \geq -C_1\left(\int_M v^{-n}\dvol_g\right)^{\frac{2m}{n(m+n)}}\left(\int_M v^m\dvol_g\right)^{\frac{2}{m+n}} . \qedhere \]
\end{proof}

This result implies that one can meaningfully ask whether minimizers of the weighted Yamabe constant exist.  We shall not treat this question in full generality here, but only present certain partial results we need in Section~\ref{sec:existence_uniqueness}.

Another useful fact is that the sign of the first eigenvalue of the weighted conformal Laplacian agrees with the sign of the weighted Yamabe constant.  In particular, it is a conformal invariant.

\begin{prop}
\label{prop:weighted_conformal_laplacian_eigenvalue}
Let $(M^n,g,v^m\dvol)$ be a compact SMMS.  The first eigenvalue $\lambda_0(L_\phi^m)$ of the weighted conformal Laplacian is positive (resp.\ zero, negative) if and only if the weighted Yamabe constant $\sigma_{1,2}(g,v^m\dvol)$ is positive (resp.\ zero, negative).
\end{prop}

\begin{proof}

This follows immediately from~\eqref{eqn:weighted_yamabe_conformal_functional} together with the standard fact that minimizers of
\[ \lambda_0(L_\phi^m) = \inf\left\{\frac{(L_\phi^mw,w)}{\lV w\rV_2^2} \colon 0\not=w\in L_1^2(M,v^m\dvol) \right\} \]
always exist and can be taken be be smooth positive functions.
\end{proof}

Our final observation about the weighted Yamabe constant is a relationship between it and the Yamabe constant, generalizing an observation of Akutagawa, Ishida and LeBrun~\cite{Akutagawa2007}.

\begin{prop}
\label{prop:generalize_akutagawa}
Let $(M^n,g,v^m\dvol)$ be a compact SMMS with positive (resp.\ nonnegative) weighted Yamabe constant.  Then the Yamabe constant $\sigma_1(g)$ of the Riemannian manifold $(M^n,g)$ is positive (resp.\ nonnegative).
\end{prop}

\begin{proof}

Since the weighted Yamabe constant and the Yamabe constant are both conformally invariant, it is enough to compare $\sigma_{1,2}(v^{-2}g,1^m\dvol_{v^{-2}g})$ and $\sigma_1(v^{-2}g)$.  Moreover, by Proposition~\ref{prop:weighted_conformal_laplacian_eigenvalue}, it suffices to compare the signs of the first eigenvalues $\lambda_0(L_\phi^m)$ and $\lambda_0(L)$, where $L$ is the conformal Laplacian.  To that end, observe that
\begin{align*}
(L_\phi^mw,w) & = \int_M \left(\frac{4(m+n-1)}{m+n-2}\lv\nabla w\rv^2 + Rw^2\right)\dvol_{v^{-2}g} \\
& \leq \int_M \left(\frac{4(n-1)}{n-2}\lv\nabla w\rv^2 + Rw^2\right)\dvol_{v^{-2}g} \\
& = (Lw,w) .
\end{align*}
The conclusion immediately follows.
\end{proof}
\subsection{The Energy of a SMMS}
\label{sec:defn_energy}

Geometrically speaking, the main motivation for the introduction of the energy is that there are no nontrivial compact quasi-Einstein SMMS with \charconstant\ zero; in other words, such a SMMS will never be a critical point of the $(m,0)$-energy functional.  For this reason, we introduce the energy as the natural conformal invariant associated to the energy functional of a SMMS with positive \charconstant; in particular, unless otherwise specified, we shall always assume that the \charconstant, when specified, is positive.

\begin{defn}
\label{defn:menergy_mu}
Let $(M^n,g,v^m\dvol)$ be a compact SMMS with \charconstant\ $\mu$ and $m<\infty$.  The \emph{$(m,\mu)$-energy $\lambda(g,v^m\dvol,\mu)$} is defined by
\[ \lambda(g,v^m\dvol,\mu) = \inf\left\{ \mW_\mu^m\left(u^{-2}g,u^{-m-n}v^m\dvol_g\right) \colon \int u^{-m-n}v^m\dvol_g=1 \right\}, \]
where the infimum is taken over all positive $u\in C^\infty(M)$.
\end{defn}

As an invariant for SMMS with \charconstant\ $\mu$, the $(m,\mu)$-energy is clearly conformally invariant.  However, it suffers a critical problem when attempting to carry out rescaling arguments \emph{while holding $v$ and $\mu$ fixed}: It is not scale invariant in this sense.  To overcome this, we define the $m$-energy by minimizing over all positive multiples of $v$, or equivalently, all positive multiples of $\mu$.

\begin{defn}
\label{defn:menergy}
Let $(M^n,g,v^m\dvol)$ be a compact SMMS with $m<\infty$.  The \emph{$m$-energy $\lambda(g,v^m\dvol)$} is defined by
\[ \lambda(g,v^m\dvol) = \inf_{\tau>0} \lambda(g,(\sqrt\tau v)^m\dvol_g,1) . \]
\end{defn}

Note that, by the conformal invariance of the $(m,\mu)$-energy,
\[ \lambda\left(cg,(\sqrt{c\tau}v)^m\dvol_g,1\right) = \lambda\left(g,(\sqrt\tau v)^m\dvol_g,1\right), \]
so the parameter $\tau$ scales relative to the metric in the same way as it does in Perelman's definitions; see below.

As stated, our definition of the $m$-energy does not make sense when $m=\infty$.  To correct this, analogous to Proposition~\ref{prop:energy_functional_smooth}, we will introduce a renormalization of the $m$-energy.  For the moment, we simply \emph{define} the $\infty$-energy as follows.

\begin{defn}
\label{defn:menergy_infty}
Let $(M^n,g,e^{-\phi}\dvol,\infty)$ be a compact SMMS.  The \emph{$\infty$-energy $\lambda(g,e^{-\phi}\dvol,\infty)$} is defined by
\[ \lambda(g,e^{-\phi}\dvol,\infty) = \inf_{\substack{f\in C^\infty(M)\\\tau>0}} \left\{ \mW_{\tau^{-1}}^m\left(g,\tau^{-\frac{n-2}{2}}e^{-f-\phi}\right) \colon \int \tau^{-\frac{n}{2}}e^{-f-\phi}\dvol=1\right\} \]
\end{defn}

As is well-known in the context of finding minimizers of Gagliardo-Nirenberg-Sobolev inequalities (cf.\ \cite{DelPinoDolbeault2002}), our definition of the $m$-energy for $m<\infty$ is essentially equivalent to the weighted Yamabe constant.

\begin{prop}
\label{prop:yamabe_and_energy}
Let $(M^n,g,v^m\dvol)$ be a compact SMMS with $m<\infty$.  Then
\begin{equation}
\label{eqn:lambda_sigma_inequality}
\lambda(g,v^m\dvol) = \begin{cases}
                        -\infty & \mbox{if } \sigma_{1,2}(g,v^m\dvol)<0 \\
                        0 & \mbox{if } \sigma_{1,2}(g,v^m\dvol)=0 \\
                        (m+n)\left(\frac{\sigma_{1,2}(g,v^m\dvol)}{n}\right)^{\frac{n}{m+n}} & \mbox{if } \sigma_{1,2}(g,v^m\dvol)>0 .
                        \end{cases}
\end{equation}
\end{prop}

\begin{remark}

When $m=\infty$, it is no longer true that one has an explicit formula relating the $m$-energy and the weighted Yamabe constant.  Nevertheless, Perelman showed~\cite{Perelman1} that the $m$-energy is finite (i.e.\ $\lambda(g,e^{-\phi}\dvol,m)\not=-\infty$) if and only if the weighted Yamabe constant is positive.
\end{remark}

\begin{proof}

From Definition~\ref{defn:menergy_mu} it is clear that the functional defining the $m$-energy is homogeneous in $u$.  Writing $w=u^{-\frac{m+n-2}{2}}$ and denoting $\lambda=\lambda(g,v^m\dvol)$, it then follows that
\begin{equation}
\label{eqn:menergy_functional}
\lambda = \inf_{\tau>0}\left\{ \frac{\tau^{\frac{m}{m+n}}(L_\phi^mw,w) + m\tau^{-\frac{n}{m+n}}\lV wv^{-1}\rV_2^2}{\lV w\rV_p^2} \colon 0\not=w\in L_1^2(M,v^m\dvol) \right\} .
\end{equation}
Denoting the weighted Yamabe constant by $\sigma=\sigma_{1,2}(g,v^m\dvol)$, it follows from~\eqref{eqn:weighted_yamabe_conformal_functional} that for any $w\in L_1^2(M,v^m\dvol)$ with $\lV w\rV_p=1$,
\begin{equation}
\label{eqn:first_lower_bound}
\lambda \geq \sigma\tau^{\frac{m}{m+n}}\lV wv^{-1}\rV_2^{-\frac{2m}{n}} + m\tau^{-\frac{n}{m+n}}\lV wv^{-1}\rV_2^2 .
\end{equation}
In particular, if $\sigma\geq0$, minimizing~\eqref{eqn:first_lower_bound} in $\tau$ yields
\[ \lambda\geq(m+n)\left(\frac{\sigma}{n}\right)^{\frac{n}{m+n}} . \]

Conversely, letting $w\in C^\infty(M)$ be a positive function such that
\[ (L_\phi^mw,w)\lV wv^{-1}\rV_2^{2m/n} \leq \sigma + \varepsilon \]
and $\lV w\rV_p=1$, we see that
\begin{equation*}
\lambda \leq (\sigma+\varepsilon)\tau^{\frac{m}{m+n}}\lV wv^{-1}\rV_2^{-\frac{2m}{n}} + m\tau^{-\frac{n}{m+n}}\lV wv^{-1}\rv_2^2 .
\end{equation*}
Minimizing again in $\tau$ and then taking the limit $\varepsilon\to0$ then yields the desired result.
\end{proof}

This result can be seen as a partial motivation for the following definition of the renormalized $m$-energy of a SMMS.

\begin{defn}
\label{defn:renormalized_energy}
Let $(M^n,g,v^m\dvol)$ be a compact SMMS with $m\geq 0$.  We define the \emph{renormalized $m$-energy $\olambda(g,v^m\dvol)$} as follows:
\begin{enumerate}
\item If $\sigma_{1,2}(g,v^m\dvol)\leq 0$, then
\[ \olambda_m(g,v^m\dvol) = 0 . \]
\item If $\sigma_{1,2}(g,v^m\dvol)>0$ and $m<\infty$, then
\[ \olambda_m(g,v^m\dvol) = (2\pi e)^{-\frac{n}{2}}\left(\frac{\lambda(g,v^m\dvol)}{m+n}\right)^{\frac{m+n}{2}} . \]
\item If $\sigma_{1,2}(g,v^m\dvol)>0$ and $m=\infty$, then
\[ \olambda(g,e^{-\phi}\dvol,\infty) = (2\pi)^{-\frac{n}{2}}\exp\left(\frac{\lambda(g,e^{-\phi}\dvol,\infty)}{2}\right) . \]
\end{enumerate}
\end{defn}

\begin{remark}
Proposition~\ref{prop:yamabe_and_energy} implies that if $(M^n,g,v^m\dvol)$ is a compact SMMS with $m<\infty$ and positive weighted Yamabe constant, then the renormalized $m$-energy satisfies
\[ \olambda(g,v^m\dvol) = \left(\frac{\sigma_{1,2}(g,v^m\dvol)}{2\pi ne}\right)^{n/2} . \]
\end{remark}

\begin{remark}
Definition~\ref{defn:renormalized_energy} is such that if $(M^n,g,e^{-\phi}\dvol,\infty)$ is a compact gradient Ricci soliton, then the renormalized energy $\olambda$ agrees with the central density $\Theta$ introduced by Cao, Hamilton and Ilmanen~\cite{CaoHamiltonIlmanen}.  To see this, recall that the $\infty$-energy and Perelman's $\nu$ entropy are defined by
\begin{align*}
\lambda(g,\dvol) & = \inf \left\{ \int \left[ \tau(R+\lv\nabla f\rv^2) + 2(f-n)\right]\tau^{-\frac{n}{2}}e^{-f}\dvol \colon \int \tau^{-\frac{n}{2}}e^{-f}\dvol = 1 \right\} \\
\nu(g) & = \inf \left\{ \int \left[ \tau(R+\lv\nabla f\rv^2) + f-n\right](4\pi\tau)^{-\frac{n}{2}}e^{-f}\dvol\colon \int (4\pi\tau)^{-\frac{n}{2}}e^{-f}\dvol=1 \right\} .
\end{align*}
It then follows easily that
\[ \nu(g) = \frac{1}{2}\lambda(g,\dvol) - \frac{n}{2}\log(2\pi) . \]
On the other hand, Cao, Hamilton and Ilmanen observed that for a gradient Ricci soliton, $\Theta=e^{\nu(g)}$, whence follows
\[ \Theta = (2\pi)^{-\frac{n}{2}}\exp\left(\frac{\lambda(g,\dvol)}{2}\right), \]
as claimed.
\end{remark}

Again, the main purpose of introducing the renormalized energy is to have a definition which is continuous for $m\in[0,\infty]$.

\begin{prop}
Let $(M^n,g)$ be a Riemannian manifold and fix a positive function $v\in C^\infty(M)$.  Then the map
\[ m \mapsto \olambda(g,v^m\dvol) \]
is continuous for $m\in[0,\infty]$.
\end{prop}

\begin{proof}

Using the conformal invariance of the $m$ and the $(m,\mu)$-energy, we may assume without loss of generality that $v\equiv 1$.  Analogous to~\eqref{eqn:menergy_functional}, for any fixed $\tau>0$ it holds that
\[ \lambda\left(g,(\sqrt\tau)^m\dvol,1\right) = \inf\left\{ \int \left[wL_\phi^m w + m\tau^{-1}w^2\right]\tau^{\frac{m}{2}} \colon \int w^{\frac{2(m+n)}{m+n-2}}\tau^{\frac{m}{2}}\dvol = 1 \right\} , \]
where the infimum is taken over all positive $w\in C^\infty(M)$ satisfying the specified constraint.  In particular, by absorbing a factor of $\tau^{\frac{m+n}{2}}$ into the variable $w$, we may equivalently write
\begin{equation}
\label{eqn:menergy_perelman}
\lambda\left(g,(\sqrt\tau)^m\dvol,1\right) = \inf_{w\in\mC(\tau)}\left\{ \int \left[\tau wL_\phi^m w + mw^2\right]\tau^{-\frac{n}{2}}\dvol \right\}
\end{equation}
for
\[ \mC(\tau) = \left\{ 0<w\in C^\infty(M) \colon \int w^{\frac{2(m+n)}{m+n-2}}\tau^{-\frac{n}{2}}\dvol = 1 \right\} . \]
Define
\[ \tilde\lambda\left(g,(\sqrt\tau)^m\dvol\right) = \inf_{w\in\mC(\tau)}\left\{ \int \left[\tau wL_\phi^m w + mw^2 - (m+2n)w^{\frac{2(m+n)}{m+n-2}}\right] \tau^{-\frac{n}{2}}\dvol \right\} . \]
On the one hand, it is clear that
\[ \lambda\left(g,(\sqrt\tau)^m\dvol,1\right) = \tilde\lambda\left(g,(\sqrt\tau)^m\dvol\right) + m + 2n . \]
On the other hand, analogous to Proposition~\ref{prop:energy_functional_smooth}, it is easy to see that
\[ \lim_{m\to\infty} \tilde\lambda\left(g,(\sqrt\tau)^m\dvol\right) = \inf_{w\in\mC(\tau)} \mW_{\mu^{-1}}^\infty\left(g,\tau^{-\frac{n-2}{2}}w\right) . \]
In particular, we see that
\begin{align*}
\lim_{m\to\infty}\olambda(g,1^m\dvol) & = \lim_{m\to\infty} (2\pi e)^{-\frac{n}{2}}\left(\frac{\inf_\tau \lambda(g,(\sqrt\tau)^m,1\dvol)}{m+n}\right)^{\frac{m+n}{2}} \\
& = \lim_{m\to\infty} (2\pi e)^{-\frac{n}{2}}\left(\frac{\inf_\tau \tilde\lambda(g,(\sqrt\tau)^m\dvol)+m+2n}{m+n}\right)^{\frac{m+n}{2}} \\
& = (2\pi)^{-\frac{n}{2}}\exp\left(\frac{\lambda(g,\dvol,\infty)}{2}\right),
\end{align*}
as desired.
\end{proof}

Finally, either using Proposition~\ref{prop:qe_var} or a straightforward direct computation, we have the following variational formulae for the $(m,1)$-energy functional restricted to volume-preserving conformal variations, as needed to understand minimizers of the $m$-energy.

\begin{prop}
\label{prop:sc_var}
Let $(M^n,g,v^m\dvol)$ be a compact SMMS and fix $\tau>0$ and a positive function $w\in C^\infty(M)$ such that $w\in\mC(\tau)$.  Let $(w_s,\tau_s)$ be a variation of $(w,\tau)$ such that $w_s\in\mC(\tau_s)$, and define
\[ \dot w = \frac{d}{ds}w_s\rv_{s=0}, \quad c = \frac{d}{ds}\left(\log\tau^{\frac{m}{2}}\right)\rv_{s=0} . \]
Denoting $u=w^{-\frac{2}{m+n-2}}$ and
\[ \delta\mW_1^m := \frac{d}{ds}\mW_1^m\left(u_s^{-2}g,u_s^{-m-n}(\sqrt\tau_s)^mv^m\dvol_g\right)\rv_{s=0} , \]
it holds that
\begin{equation}
\label{eqn:sc_var}
\begin{split}
\delta\mW_1^m & = 2\int_M \left(L_\phi^m w + m\tau^{-1}v^{-2}w\right)\dot w\,\tau^{\frac{m}{2}}v^m\dvol \\
& \quad + \int_M \left(wL_\phi^m w + (m-2)\tau^{-1}v^{-2}w^2\right)c\,\tau^{\frac{m}{2}}v^m\dvol_g .
\end{split}
\end{equation}
and
\begin{equation}
\label{eqn:sc_var_meas}
0 = c+ \frac{2(m+n)}{m+n-2}\int_M \dot{w} w^{\frac{m+n+2}{m+n-2}}\tau^{\frac{m}{2}}v^m\dvol .
\end{equation}
\end{prop}

Using Corollary~\ref{cor:weighted_conformal_laplacian_scalar}, we have the following corollaries of Proposition~\ref{prop:sc_var}.

\begin{cor}
\label{cor:sc_crit_mmu}
Let $(M^n,g,v^m\dvol)$ be a compact SMMS and suppose that $1\in\mC(1)$ is a critical point of the $(m,1)$-energy; i.e.\ $\mW_1^m(g,v^m\dvol) = \lambda(g,v^m\dvol,1)$.  Then
\begin{equation}
\label{eqn:sc_crit_scalar_mmu} R_\phi^m + mv^{-2} = \lambda(g,v^m\dvol) .
\end{equation}
\end{cor}

\begin{cor}
\label{cor:sc_crit}
Let $(M^n,g,v^m\dvol)$ be a compact SMMS with $m>0$ and suppose that $1\in\mC(1)$ is a critical point of the $m$-energy; i.e.\ $\mW_1^m(g,v^m\dvol) = \lambda(g,v^m\dvol)$.  Then
\begin{align}
\label{eqn:sc_crit_scalar} R_\phi^m + mv^{-2} & = \lambda(g,v^m\dvol) \\
\label{eqn:sc_crit_laplacian} \int_M v^{m-2}\dvol_g & = \frac{\lambda(g,v^m\dvol)}{m+n} .
\end{align}
\end{cor}                   % Energy Definition
\section{Towards the Existence and Uniqueness of Minimizers}
\label{sec:existence_uniqueness}

In this section, we consider some basic questions about the existence and the uniqueness of minimizers of the $(m,\mu)$-energy as they will be needed to establish our precompactness theorem.  We expect that one can extend our results to the $m$-energy itself, but we do not pursue this question here.  Also, for convenience and unless otherwise specified, we shall henceforth assume that $m>0$.

First, let us consider the existence question for the $(m,\mu)$-energy, which is anyway a necessary first step in establishing the existence of minimizers of the weighted Yamabe constant in general (cf.\ \cite{DelPinoDolbeault2002,Perelman1}).

\begin{lem}
\label{lem:mphi_energy_minimizer}
Let $(M^n,g,v^m\dvol)$ be a SMMS with \charconstant\ $\mu>0$.  Then there exists a positive function $u\in C^\infty(M)$ such that
\begin{align}
\label{eqn:menergy_minimizer} \mW_\mu^m(u^{-2}g,u^{-m-n}v^m\dvol_g) & = \lambda(g,v^m\dvol,\mu) \\
\label{eqn:menergy_normalization} \int_M u^{-m-n}v^m\dvol_g & = 1 
\end{align}
\end{lem}

\begin{proof}

Set $w=u^{-\frac{m+n-2}{2}}$, so that fixing $\tau=1$ in~\eqref{eqn:menergy_functional} implies that
\[ \mW_\mu^m(u^{-2}g,u^{-m-n}v^m\dvol_g) \geq \frac{4(m+n-1)}{m+n-2}\lV\nabla w\rV_2^2 - C_1\lV w\rV_2^2 \geq -C_2, \]
where the constants $C_1,C_2>0$ depend only on $(M^n,g,v^m\dvol)$ and $\mu$.  We may thus choose a minimizing sequence $\{w_i\}\subset C^\infty(M)$ of positive functions such that $\lV w\rV_p=1$ for $p=\frac{2(m+n)}{m+n-2}$.  Since $m>0$, the embedding $L_1^2(M,v^m\dvol)\hookrightarrow L^p(M,v^m\dvol)$ is compact, and hence, after taking a subsequence if necessary, $w_i$ converges strongly to some nonnegative $w\in L^p(M,v^m\dvol)$ which realizes the equality~\eqref{eqn:menergy_minimizer}.  By Proposition~\ref{prop:sc_var}, we know that $w$ is a weak solution of
\[ -\frac{4(m+n-1)}{m+n-2}\Delta_\phi w + (R_\phi^m+m\mu v^{-2})w = \lambda w^{\frac{m+n+2}{m+n-2}} , \]
where we have written $\lambda=\lambda(g,v^m\dvol,\mu)$.  Standard arguments using elliptic regularity theory then imply that $0<w\in C^\infty(M)$, as desired.
\end{proof}

Second, Obata's classification of the minimizers of the Yamabe constant of an Einstein manifold~\cite{Obata1971} extends to the $(m,\mu)$-energy in a natural way.

\begin{thm}
\label{thm:convexity}
Let $(M^n,g,v^m\dvol)$ be a compact quasi-Einstein SMMS with \charconstant\ $\mu$ and $\Vol_\phi(M)=1$.  Suppose that $w$ is a minimizer of the $(m,\mu)$-energy as in Lemma~\ref{lem:mphi_energy_minimizer}.  Then either $w=1$ or $m=0$ and $(M^n,g)$ is isometric to the standard sphere with $\Ric_g=g$ and $w^{\frac{4}{n-2}}g$ is an Einstein metric.
\end{thm}

\begin{proof}

When $m=0$, this is precisely Obata's result~\cite{Obata1971}, while when $m=\infty$, this is a result of Perelman~\cite{Perelman1}.  We thus assume $m\in(0,\infty)$.

The key observation is that, by a simple algebraic computation, the variational formula from Proposition~\ref{prop:gen_var} can be rewritten as
\begin{align*}
\delta\mW_\mu^m & = -\int_M \lp \Ric_\phi^m - \frac{R_\phi^m+m\mu v^{-2}}{m+n}g, h + \frac{2}{m}\psi g\rp v^m\dvol_g \\
& \quad - \frac{m+n-2}{m+n}\int_M \left(R_\phi^m+m\mu v^{-2}\right)\left(\psi-\frac{1}{2}\tr h\right) v^m\dvol_g .
\end{align*}
In particular, since the $(m,\mu)$-energy functional is invariant under diffeomorphisms, this vanishes when we take $h=L_Xg$ and $\psi=X\phi$; i.e.\ the general Pohozaev-type formula
\begin{equation}
\label{eqn:pohozaev}
\begin{split}
0 & = \frac{m+n-2}{m+n}\int_M X\left(R_\phi^m+m\mu v^{-2}\right) v^m\dvol_g \\
& \quad + \int_M \lp\Ric_\phi^m - \frac{R_\phi^m+m\mu v^{-2}}{m+n}g, v^2L_X(v^{-2}g)\rp v^m\dvol_g
\end{split}
\end{equation}
holds by a straightforward integration-by-parts argument.

To prove the theorem, we apply~\eqref{eqn:pohozaev} to the SMMS
\begin{equation}
\label{eqn:smms_w}
\left( M^n, w^{\frac{4}{m+n-2}}g, w^{\frac{2(m+n)}{m+n-2}}v^m\dvol_g \right) .
\end{equation}
The assumption that $w$ is a minimizer of the $(m,1)$-energy implies, using Corollary~\ref{cor:sc_crit_mmu}, that $\hat R_\phi^m+m\mu\hat v^{-2}$ is constant, where the hats mean that we compute using~\eqref{eqn:smms_w} rather than $(M^n,g,v^m\dvol)$.  For convenience, define
\[ E_\phi^m := \Ric_\phi^m - \frac{R_\phi^m+m\mu v^{-2}}{m+n}g , \]
and let $\hat E_\phi^m$ denote the same quantity computed for~\eqref{eqn:smms_w}.  Using Proposition~\ref{prop:curvature_conformal}, we see that
\[ \hat E_\phi^m = E_\phi^m + (m+n-2)u^{-1}\left(\nabla^2 u - \frac{1}{m+n}\Delta_\phi u\,g\right) . \]
On the other hand, the assumption that $(M^n,g,v^m\dvol)$ is quasi-Einstein with \charconstant\ $\mu$ implies that $E_\phi^m=0$, whence~\eqref{eqn:pohozaev} implies that
\begin{equation}
\label{eqn:pohozaev_obata}
0 = (m+n-2)\int_M \lp \nabla^2 u - \frac{1}{m+n}\Delta_\phi u\,g, L_Xg + \frac{2}{m}X\phi\,g\rp u^{1-m-n}v^m\dvol_g, 
\end{equation}
where we are now computing with respect to $(M^n,g,v^m\dvol)$.

Finally, observe that if we set $X=\frac{1}{2}\nabla u$, it holds that
\begin{align*}
A & := L_Xg + \frac{2}{m}X\phi\,g \\
& = \nabla^2 u + \frac{1}{m}\lp\nabla u,\nabla\phi\rp\,g \\
& = \left(\nabla^2 u - \frac{1}{m+n}\Delta_\phi u\,g\right) + \frac{1}{m}\tr_g\left(\nabla^2 u - \frac{1}{m+n}\Delta_\phi u\,g\right) g .
\end{align*}
It thus follows from~\eqref{eqn:pohozaev_obata} that $\nabla^2u-\frac{1}{m+n}\Delta_\phi u\,g=0$.  Taking the trace and applying the maximum principle then implies that $u$ is constant, as desired.
\end{proof}

\begin{remark}
The above proof in fact yields the cases $m=0$ and $m=\infty$ without change.  However, when $m=0$, one must be careful that $E_\phi^m$ is always tracefree, whence~\eqref{eqn:pohozaev} only allows one to conclude that the tracefree part of $\nabla^2u$ vanishes.
\end{remark}

Finally, while we cannot establish that minimizers of the $m$-energy on a quasi-Einstein SMMS are unique, and thereby compute their $m$-energy, we can nevertheless establish lower bounds for the quasi-Einstein constant, the \charconstant, and the weighted volume of a quasi-Einstein SMMS in terms of its $m$-energy.

\begin{prop}
\label{prop:compute_energy_quasi_einstein}
Let $(M^n,g,v^m\dvol)$ be a compact quasi-Einstein SMMS with quasi-Einstein constant $\lambda$ and \charconstant\ $\mu>0$.  Then
\[ \olambda(g,v^m\dvol) \leq (2\pi e)^{-\frac{n}{2}} \lambda^{\frac{m+n}{2}}\mu^{-\frac{m}{2}} \Vol_\phi(M) . \]
\end{prop}

\begin{proof}

By~\eqref{eqn:menergy_functional} applied with $w=1$ and $\tau=\mu^{-1}$, it follows that
\[ \lambda(g,v^m\dvol) \leq \frac{\mu^{-\frac{m}{m+n}}\int R_\phi^mv^m\dvol + m\mu^{\frac{n}{m+n}}\int v^{m-2}\dvol}{\Vol_\phi(M)^{2/p}} . \]
On the other hand, \eqref{eqn:bianchi} implies that
\[ \int R_\phi^mv^m\dvol = (m+n)\lambda\Vol_\phi(M) - m\mu\int v^{m-2}\dvol . \]
It thus follows that
\[ \lambda(g,v^m\dvol) \leq (m+n)\lambda\mu^{-\frac{m}{m+n}}\Vol_\phi(M)^{\frac{2}{m+n}} . \]
The result then follows from the definition of the renormalized energy .
\end{proof}   % Existence/Uniqueness
\section{Relating the Energy to the Geometry of a SMMS}
\label{sec:kappa_noncollapsing}

An important property of the Yamabe constant and Perelman's $\nu$-entropy is their relationship to Sobolev and logarithmic Sobolev inequalities.  As far as the author can tell, sharp weighted Sobolev inequalities on SMMS are not yet completely understood in general.  However, in the presence of lower bounds on the Bakry-\'Emery Ricci tensor, the following result is known (cf.\ \cite{Ledoux2000,Villani2009}):

\begin{thm}
\label{thm:sharp_sobolev}
Let $(M^n,g,v^m\dvol)$ be a compact SMMS such that $\Vol_\phi(M)=1$ and $\Ric_\phi^m\geq Kg>0$.  Then for all $w\in L_1^2(M,v^m\dvol)$,
\begin{equation}
\label{eqn:sharp_sobolev}
\lV w\rV_p^2 \leq \frac{4(m+n-1)}{(m+n)(m+n-2)K} \lV\nabla w\rV_2^2 + \lV w\rV_2^2,
\end{equation}
where all norms are computed with respect to $v^m\dvol$.
\end{thm}

On the other hand, a similar inequality holds for SMMS which already minimize the energy.

\begin{prop}
Let $(M^n,g,v^m\dvol)$ be a compact SMMS with positive $m$-energy $\lambda$, and suppose that $(w,\tau)=(1,1)$ minimizes~\eqref{eqn:menergy_functional}.  Then for all $w\in L_1^2(M,v^m\dvol)$,
\[ \lV w\rV_p^2 \leq \frac{4(m+n-1)}{(m+n-2)\lambda}\lV\nabla w\rV_2^2 + \lV w\rV_2^2 . \]
\end{prop}

\begin{proof}

From~\eqref{eqn:menergy_functional} with $\tau=1$, it follows that
\[ \lambda\lV w\rV_p^2 \leq \frac{4(m+n-1)}{m+n-2}\lV\nabla w\rV_2^2 + \int \left( R_\phi^m+mv^{-2}\right)w^2\,v^m\dvol . \]
On the other hand, by Corollary~\ref{cor:sc_crit}, the assumption that $(1,1)$ minimizes~\eqref{eqn:menergy_functional} implies that
\[ R_\phi^m+mv^{-2}=\lambda . \]
Combining these equations yields the result.
\end{proof}

Another important feature of the relationship between the Yamabe constant and Perelman's $\nu$-entropy to Sobolev-type inequalities is the relationship between the latter and the isoperimetric inequality.  Geometrically, this manifests itself in the $\kappa$-noncollapsing property of manifolds with positive Yamabe constant (cf.\ \cite{Anderson1989}) and positive $\nu$-entropy~\cite{Perelman1}.  Unsurprisingly, similar results hold true for SMMS with positive $m$-energy, as we illustrate in two ways.

First, we have the following global estimate for the volume in terms of the renormalized $m$-energy, which is valid given the existence of a minimizer of the $m$-energy.

\begin{prop}
\label{prop:global_volume_bound}
Let $(M^n,g,1^m\dvol)$ be a compact SMMS with positive renormalized energy $\olambda$.  Suppose additionally that $(w,\tau)$ is a minimizer of the $m$-energy exists.  Then
\begin{equation}
\label{eqn:energy_volume}
\Vol(M) \geq (2\pi e\tau)^{n/2} \olambda .
\end{equation}
\end{prop}

\begin{proof}

Set $u=w^{-\frac{2}{m+n-2}}$, so that the assumptions imply that
\[ \left( M^n, u^{-2}g, \left(\sqrt{\tau}u^{-1}\right)^m\dvol_{u^{-2}g} \right) \]
is such that the constant function $1$ is a minimizer of the $m$-energy.  It follows from Corollary~\ref{cor:sc_crit}, H\"older's inequality, and the normalization of $w$ that
\[ \frac{\lambda(g,1^m\dvol)}{m+n} = \int \tau^{\frac{m-2}{2}}u^{2-m-n}\dvol_g \leq \tau^{-\frac{n}{m+n}}\Vol(M)^{\frac{2}{m+n}} . \]
The result then follows from the definition of the renormalized energy.
\end{proof}

While this result is not immediately applicable here, it is nevertheless the inspiration for the following global volume estimate for conformally quasi-Einstein SMMS (cf.\ Proposition~\ref{prop:compute_energy_quasi_einstein}, Proposition~\ref{prop:global_volume_bound}).

\begin{prop}
\label{prop:global_volume_bound_qe}
Let $(M^n,g,1^m\dvol)$ be a compact SMMS with \charconstant\ one and suppose that $u\in C^\infty(M)$ is a positive quasi-Einstein scale with quasi-Einstein constant $\lambda$, normalized to have unit volume.  Then
\[ \Vol(M)^{\frac{2}{m+n}} \geq \lambda . \]
\end{prop}

\begin{proof}

Subtracting~\eqref{eqn:qe_mu_uv} from~\eqref{eqn:qe_lambda_uv} and integrating with respect to $u^{-m-n}\dvol_g$, it follows that
\[ \lambda = \int u^{2-m-n}\dvol_g . \]
The result then follows from H\"older's inequality as above.
\end{proof}

Second, we can establish the $\kappa$-noncollapsing property for compact quasi-Einstein SMMS with positive $m$-energy (cf.\ \cite{Perelman1}).  First, let us recall the definition of a $\kappa$-noncollapsed manifold.

\begin{defn}
Given $\kappa,d>0$, a Riemannian manifold $(M,g)$ is said to be \emph{$\kappa$-noncollapsed below the scale} $d$ if, whenever $B=B(x,r)$ is such that $r<d$ and $R\leq r^{-2}$ on $B$, then $\Vol(B)\geq \kappa r^{n}$.
\end{defn}

By way of motivation, let us first establish a $\kappa$-noncollapsing result assuming the existence of a minimizer, but without the quasi-Einstein assumption.

\begin{thm}
\label{thm:kappa}
Let $(M^n,g,1^m\dvol)$ be a compact SMMS such that the renormalized $m$-energy $\olambda\geq C>0$, and suppose additionally that there exists a positive minimizer $(w,\tau)$ of~\eqref{eqn:menergy_functional}.  Then for every $d>0$, there is a constant $\kappa(C,d,n)>0$ such that $(M,g)$ is $\kappa$-noncollapsed below the scale $d\tau^{1/2}$.
\end{thm}

\begin{proof}

Suppose the theorem is false, and so there is a sequence $\{B_k=B(x_k,r_k)\}$ inside the manifolds $(M_k,g_k)$ such that $R\leq r_k^{-2}$ in $B_k$ but $r_k^{-n}\Vol(B_k)\to 0$ as $k\to\infty$, and that additionally, there are $0\leq m_k\leq \infty$ such that $\olambda_{m_k}(g_k,1^{m_k}\dvol_{g_k})\geq C$.  As in~\cite{KleinerLott2008}, by redefining $r_k=r_k/2$ if necessary, we may assume that
\begin{equation}
\label{eqn:kappa_doubling}
\frac{\Vol(B_k)}{\Vol(B(x_k,r_k/2))} < 3^n .
\end{equation}

Suppose first that $\partial B_k=\emptyset$ for all $k$.  Then $B_k=M_k$, and we may apply Proposition~\ref{prop:global_volume_bound} to see that $r_k^{-n}\Vol(B_k) \geq d^{-n} C$, a contradiction.

Thus, taking a subsequence if necessary, we may assume that $\partial B_k\not=\emptyset$ for all $k$.  For simplicity, suppose that this sequence is such that $m_k<\infty$ for all $k$, and define the test functions
\[ \xi_k=\exp\left(-\frac{m+n-2}{2(m+n)}c_k\right)\eta(d(x_k,\cdot)/r_k) , \]
where $\eta\colon[0,\infty)\to[0,1]$ is a smooth cutoff function such that $\eta\equiv 1$ on $[0,1/2]$ and $\eta\equiv 0$ on $[1,\infty)$, and $c_k$ is chosen such that
\begin{equation}
\label{eqn:kappa_normalization}
\int_M r_k^{-n} \xi_k^{2\frac{m+n}{m+n-2}} \dvol = 1 ,
\end{equation}
where we suppress the dependence of $m$ on $k$ in the interest of readability.

This normalization together with the choice of cutoff function has three important consequences.  First, we have that
\begin{equation}
\label{eqn:kappa_eqn1}
e^{c_k} = \int_M r_k^{-n} \eta^{2\frac{m+n}{m+n-2}} \dvol \leq r_k^{-n} \Vol(B_k) .
\end{equation}
Second, using H\"older's inequality, we have that
\begin{equation}
\label{eqn:kappa_eqn2}
\int_M \xi_k^2 r_k^{-n} \dvol \leq \left( r_k^{-n} \Vol(B_k) \right)^{\frac{2}{m+n}} .
\end{equation}
Third, using the doubling assumption~\eqref{eqn:kappa_doubling}, we have that
\begin{equation}
\label{eqn:kappa_eqn3}
e^{-c_k}r_k^{-n}\Vol(B_k) < 3^n e^{-c_k} r_k^{-n}\Vol(B(x_k,r_k/2)) \leq 3^n .
\end{equation}
Since $\eta$ is fixed, using $\tau=r_k^2$ and $w=\xi_k$ in~\eqref{eqn:menergy_perelman}, we see that
\begin{align*}
\lambda(g,r_k^m\dvol,1) & \leq \int_M \left[4\frac{m+n-1}{m+n-2}r_k^2|\nabla \xi_k|^2 + (r_k^2 R + m)\xi_k^2\right]r_k^{-n} \dvol \\
& = m\int_M \xi_k^2 r_k^{-n} \dvol \\
& \quad + e^{-\frac{m+n-2}{m+n}c_k} \int_M \left( 4\frac{m+n-1}{m+n-2}(\eta^\prime)^2 + \eta^2\right) r_k^{-n} \dvol \\
& \leq \left( r_k^{-n} \Vol(B_k) \right)^{\frac{2}{m+n}} \left(m+c(n)\right) ,
\end{align*}
where $c(n)=3^n(4\frac{n-1}{n-2}\sup(\eta^\prime)^2+1)$ and we have used~\eqref{eqn:kappa_eqn1}, \eqref{eqn:kappa_eqn2}, and~\eqref{eqn:kappa_eqn3} in the final inequality.

Using the definition of $\olambda_m(g_k,0)$, we then see that
\[ (2\pi e)^{n/2}\olambda_m(g_k,0) \leq \left(1+\frac{c(n)}{m+n}\right)^{\frac{m+n}{2}} r_k^{-n}\Vol(B_k) \leq e^{c(n)/2} r_k^{-n} \Vol(B_k) , \]
where the second inequality follows from the elementary estimate $(1+c/m)^m\leq e^c$.  Letting $k\to0$, it then follows that $\olambda_m(g_k,0)$ tends to zero, a contradiction.
\end{proof}

The only place we used the existence of a minimizer was to rule out the possibility that the sequence of manifolds was globally collapsing via Proposition~\ref{prop:global_volume_bound}.  Thus we can apply the same argument to establish a $\kappa$-noncollapsing result for compact quasi-Einstein SMMS.  In fact, we can remove the dependence of $\kappa$ on $d$ (cf.\ \cite{Perelman1}).

\begin{defn}
A Riemannian manifold $(M,g)$ is \emph{$\kappa$-noncollapsed on all scales} if it is $\kappa$-noncollapsed below the scale $d$ for all $d>0$.
\end{defn}

\begin{thm}
\label{thm:kappa_qe}
Given $C>0$, $3\leq n\in\bN$, there is a constant $\kappa(C,n)>0$ such that any compact SMMS $(M^n,g,1^m\dvol)$ with \charconstant\ one and renormalized $m$-energy $\olambda\geq C$ which admits a quasi-Einstein scale is $\kappa$-noncollapsed on all scales.
\end{thm}

\begin{proof}

Following the proof of Theorem~\ref{thm:kappa}, it suffices to show that if $M=B(x,r)$ satisfies $r^2 R\leq 1$, then $r^{-n}\Vol(M)$ is uniformly bounded below.  Indeed, by Proposition~\ref{prop:global_volume_bound_qe}, it suffices to show that $r^{-2}$ is uniformly bounded below.  This follows from~\eqref{eqn:qe_mu_uv} by noting that if $p$ maximizes $u$, then $R(p)\geq n$, whence $r^{-2}\geq n$.
\end{proof}

While Theorem~\ref{thm:kappa_qe} is local, in the presence of global scalar curvature bounds, we can also achieve global estimates on the growth rate of geodesic balls.  First, we have the following simple volume estimate for SMMS with bounded scalar curvature:

\begin{lem}
\label{lem:vol_growth}
Let $(M^n,g,1^m\dvol)$ be a compact SMMS with \charconstant\ one and renormalized $m$-energy $\olambda\geq C>0$ which admits a quasi-Einstein scale.  Suppose additionally that the scalar curvature $R\leq C_2^2$.  Then there is a constant $c(n,C_1,C_2)>0$ such that $\Vol\left(B(x,r)\right)\geq cr$ for all $x\in M$, $r<\diam(M)$.
\end{lem}

\begin{proof}

By Theorem~\ref{thm:kappa_qe}, if $r_1\leq C_2^{-1}$, then $\Vol(B(x,r_1))\geq \kappa r_1^n$.  Given $r<\diam(M)$, we can fit at least $\frac{C_2r}{2}$ disjoint balls of radius $C_2^{-1}$ into $B(x,r)$, and so
\[ \Vol(B(x,r)) \geq \sum_{i=1}^{C_2r/2} \Vol(B(x_i,C_2^{-1})) \geq \frac{C_2^{1-n}}{2}\kappa r . \qedhere \]
\end{proof}

However, this lemma is insufficient for studying quasi-Einstein SMMS, where we only know that the scalar curvature grows at most quadratically (see Proposition~\ref{prop:bounds_no_m}).  We overcome this problem by modifying an argument of Sesum and Tian~\cite{SesumTian2008} and H.-D.\ Cao and X.-P.\ Zhu (see~\cite{Cao2009_surveyb}) to establish that the volumes of geodesic balls grows at least logarithmically in the radius.

\begin{thm}
\label{thm:vol_growth_log}
Given $c_1,c_2,c_3>0$, $3\leq n\in\bN$, there is a constant $C(n,c_1,c_2,c_3)>0$ such that whenever $(M^n,g,1^m\dvol)$ is a compact SMMS with \charconstant\ one and renormalized $m$-energy $\olambda\geq c_1$ which admits a quasi-Einstein scale $u=e^{\frac{f}{m+n-2}}\in C^\infty(M)$ such that $\lv\nabla f\rv^2,R\leq c_2r^2+c_3$, where $r$ is the distance from some fixed point $p\in M$, it holds that
\[ \Vol\left(B(p,r)\right) \geq C \log\log r \]
for all $r<\diam(M)$.
\end{thm}

\begin{proof}

Our proof is modeled after the one given by Sesum and Tian~\cite{SesumTian2008}.  Since $\olambda\geq C_1$, Theorem~\ref{thm:kappa_qe} yields a $\kappa>0$ such that $(M,g)$ is $\kappa$-noncollapsed.  As in the proof of Theorem~\ref{thm:kappa_qe}, the proof proceeds by contradiction.  To that end, let $(M_i^n,g_i,1^{m_i}\dvol_{g_i})$ be a sequence of SMMS with quasi-Einstein scales $u_i>0$ as in the theorem, and suppose that $\Vol\left(B(p_i,r_i)\right)\leq\varepsilon_i\log\log r_i$ for a sequence of points $p_i\in M_i$ and radii $r_i>0$ with $\varepsilon_i\to 0$ as $i\to\infty$.  For notational simplicity, we will henceforth drop the subscript $i$.

If $\diam(M)$ is uniformly bounded, we can apply Lemma~\ref{lem:vol_growth} and achieve a stronger linear growth rate.  Otherwise, we may assume that $\diam(M)$ is sufficiently large so that the arguments of~\cite{SesumTian2008} are valid.  Define
\[ A(r_1,r_2) = \left\{ x\in M\colon r_1 < d(p,x) < r_2 \right\}, \]
$A_{k_1,k_2}=A(2^{k_1},2^{k_2})$, and $V(k_1,k_2)=\Vol(A_{k_1,k_2})$.  Lemma~\ref{lem:vol_growth} and the quadratic growth of the scalar curvature imply
\[ V(k,k+1) \geq 2^{2(1-n)k}c_3 \]
for some constant $c_3>0$ depending only on $n$ and $C_1$.  The argument of~\cite[Claim 10]{SesumTian2008} then allows us to find $k_1,k_2$, depending only on $\varepsilon$, $n$ and $C_1$ such that $k_1+2<k_2-2$,
\begin{equation}
\label{eqn:doubling}
V(k_1,k_2)\leq \varepsilon \qquad \mbox{and} \qquad V(k_1,k_2) \leq 2^{10n}V(k_1+2,k_2-2) . 
\end{equation}
We remark here that the second property can be thought of as a ``doubling'' assumption, and is crucial in our eventual use of the $m$-energy.  Also, it is in proving this fact that we use the fact that $\diam(M)$ is sufficiently large.

As a consequence of the fundamental theorem of calculus and~\eqref{eqn:doubling}, we can find $r_1\in[2^{k_1},2^{k_1+1}]$, $r_2\in[2^{k_2-1},2^{k_2}]$ such that
\[ \Vol(\partial B(p,r_i)) \leq 2V(k_1,k_2)2^{-k_i} \]
for $i=1,2$ (cf.\ \cite{SesumTian2008}).  By Corollary~\ref{cor:duality}, we have that
\[ \Ric + \nabla^2 f + \frac{1}{m+n-2}df\otimes df = g. \]
Taking the trace and integrating with respect to $\dvol$, we see that
\[ \int_{A(r_1,r_2)} R\,\dvol \leq nV(k_1,k_2) + \int_{\partial B(x,r_1)} \lv\nabla f\rv \dvol + \int_{\partial B(x,r_2)} \lv\nabla f\rv \dvol . \]
By the assumption on $|\nabla f|$ and the definition of $r_1,r_2$, we thus conclude that
\[ \int_{A(r_1,r_2)} R\,\dvol \leq c_4 V(r_1,r_2) \]
for some constant $c_4>0$ depending only on $n,\varepsilon$, and $C_1$.

The remainder of the proof is essentially identical to the proof of Theorem~\ref{thm:kappa}, but with a different family of test functions.  Let $\eta\colon\bR\to[0,1]$ be a smooth function such that $\eta\equiv 1$ on $[2^{k_1+2},2^{k_2-2}]$, $\eta\equiv 0$ outside $[r_1,r_2]$, and $|\eta^\prime|\leq 1$.  Define the test function
\[ w=e^{\frac{m+n-2}{2(m+n)}L} \eta\circ r\]
by requiring
\[ \int_M w^{\frac{2(m+n)}{m+n-2}}\dvol = 1 . \]
This in particular implies that
\[ e^{L} V(k_1,k_2) \geq 1 , \]
and hence $L\to\infty$ as $\varepsilon\to0$.  Using the doubling property~\eqref{eqn:doubling}, a computation similar to that of the proof of Theorem~\ref{thm:kappa} yields
\[ \olambda \leq c_5 e^{-L} \to 0 \]
as $\varepsilon\to 0$, contradicting our assumption on $\olambda$ (cf.\ \cite{Cao2009_surveyb,SesumTian2008}).
\end{proof}    % Noncollapsing
\section{A Precompactness Theorem}
\label{sec:precpt}

There has recently been a lot of interest in proving precompactness theorems for Riemannian manifolds under various geometric assumptions, many of which are heavily motivated by results of Anderson~\cite{Anderson1989}, Bando, Kasue and Nakajima~\cite{BKN1989}, and Tian~\cite{Tian1990} for compact (K\"ahler-)Einstein manifolds.  Perhaps the most natural generalization of these results is to the case of compact gradient Ricci solitons.  As is well-known, compact steady or expanding gradient Ricci solitons are necessarily Einstein manifolds, and so the question is only interesting when considering compact shrinking gradient Ricci solitons.  This question was first considered by Cao and Sesum~\cite{Cao_Sesum} in the K\"ahler category, and then by X.\ Zhang~\cite{Zhang2006}, Weber~\cite{Weber2008}, and Z.\ Zhang~\cite{Zhang2009} in the Riemannian category under successively more general assumptions, with the result in~\cite{Zhang2009} as general as the one given by Anderson~\cite{Anderson1989}.  In this section, we generalize these results to quasi-Einstein SMMS, achieving a result as general as in~\cite{Anderson1989,Zhang2009}.

\begin{remark}
Haslhofer and M\"uller~\cite{HaslhoferMuller2010} have very recently removed the compactness assumption from Z.\ Zhang's result, and it should be a straightforward matter to adapt their work to our setting if one is able to remove the compactness from certain results of Section~\ref{sec:existence_uniqueness} and Section~\ref{sec:kappa_noncollapsing}.
\end{remark}

In~\cite{Case2010a}, it was pointed out that a number of general constructions for quasi-Einstein metrics actually produce smooth families parameterized by $m\in(1,\infty]$.  In light of this, the following precompactness theorem can be regarded as stating that this behavior is actually typical: Under natural geometric assumptions, one can always find a convergent subsequence of a sequence of compact quasi-Einstein SMMS with $m\to\infty$.

\begin{thm}
\label{thm:precpt}
Let $(M_i^n,g_i,1^{m_i}\dvol)$, $n\geq 4$, be a sequence of compact SMMS with \charconstant\ one and quasi-Einstein scales $u_i\in C^\infty(M_i)$, and suppose additionally that
\begin{enumerate}
\item $1<m_i\leq\infty$ and $m_i\to m>1$ as $i\to\infty$,
\item $\displaystyle\int_{M_i} u_i^{-m_i-n} 1^{m_i}\dvol_{g_i} = 1$,
\item $\olambda(g_i,1^{m_i}\dvol_{g_i})\geq C_1$
\item $\displaystyle\int_{M_i} \left|A_i\right|^{n/2} 1^{m_i}\dvol_{g_i} \leq C_2$
\end{enumerate}
for constants $C_1,C_2>0$.  Then there is a subsequence which converges in the Cheeger-Gromov sense the an orbifold $(M^n,g,1^m\dvol_g)$ with \charconstant\ one, quasi-Einstein scale $u=e^{f/(m+n-2}$, and finitely many singularities.
\end{thm}

The convergence described in Theorem~\ref{thm:precpt} is as follows (cf.\ \cite{Anderson1989,Cao_Sesum,Weber2008,Zhang2006,Zhang2008}):

\begin{defn}
A sequence of SMMS $(M_i^n,g_i,1^{m_i}\dvol_{g_i})$ with \charconstant\ one and quasi-Einstein scales $u_i=e^{f_i/(m_i+n-2)}$ \emph{converges in the Cheeger-Gromov sense} to the orbifold $(M^n,g,1^m\dvol_g)$ with \charconstant\ one, quasi-Einstein scale $u=e^{f/(m+n-2)}$, and finitely many singularities $\mS=\{p_1,\dotsc,p_k\}$ if
\begin{enumerate}
\item $m_i\to m$,
\item $(M_i,g_i)\to (M,g)$ in the Gromov-Hausdorff sense,
\item for any compact subset $K\subset M\setminus \mS$, there are compact sets $K_i\subset M_i$ and diffeomorphisms $\varphi_i\colon K\to K_i$ so that $\varphi_i^\ast g_i$ and $\varphi_i^\ast f_i$ converge in the $C^\infty$ topology to $g$ and $f$, respectively, and
\item for every $p_j\in\mS$, there is a neighborhood $U_j$ of $p_i$ which is covered by a ball $B_j\subset\bR^n$, and there is a diffeomorphism $\psi_j$ of $B_j$ such that the SMMS $(B_j,\psi_j^\ast\pi_j^\ast g,1^m\dvol)$ with \charconstant\ one admits $u=e^{\psi_j^\ast\pi_j^\ast f/m+n-2}$ as a quasi-Einstein scale, where $\pi_j\colon B_j\to U_j$ is the covering map.
\end{enumerate}
\end{defn}

Given a point $p_i\in\mS$, we will call the group $\Gamma_i$ which is such that $U_i$ is homeomorphic to $B_i/\Gamma_i$ the \emph{orbifold group}.

Before we begin proving this theorem, let us make some comments on the assumptions and our conclusion:

1) In order to prove smooth convergence, it is easier to work with a conformally quasi-Einstein SMMS using the standard measure, rather than the quasi-Einstein SMMS directly.  As mentioned in Section~\ref{sec:results}, this is because the positive elliptic Gaussian is complete when viewed as in~\eqref{eqn:hyperbolic_peg}, which is the same perspective taken in Theorem~\ref{thm:precpt}.  In particular, this will make it easier to apply the \emph{a priori} estimates from Section~\ref{sec:results}.

2) Fixing the \charconstant\ and the total volume of the conformally equivalent SMMS determined by $u_i$ to both be one removes the freedom to rescale $g\mapsto c^2g$ and $u\mapsto k^2u$, which are the two trivial sources of noncompactness in this problem (we have already fixed the size of the potential $v$ by requiring $v=1$).

3) The assumption $m,m_i>1$ is necessary in order to apply the \emph{a priori} estimates from Corollary~\ref{cor:dil}.  We do not know whether or not it is possible to weaken this assumption on $m_i$.  However, our assumptions are not enough in the case $m_i=0$ for all $i$, as our assumptions are only that we have a family of Riemannian manifolds with conformal factors $u_i$ such that the metrics $u_i^{-2}g_i$ are unit-volume Einstein metrics on $M_i$.

4) The phrasing of Theorem~\ref{thm:precpt} in terms of the $m$-energy is partly meant to emphasize its role in the underlying analysis, which is primarily based on $\kappa$-noncollapsing and its consequences.  This is the same approach taken in the work of Weber~\cite{Weber2008} and of Z.\ Zhang~\cite{Zhang2009} on limits of gradient Ricci solitons.  As we shall see in Proposition~\ref{prop:diameter_bound}, a positive lower bound for the renormalized $m$-energy of a quasi-Einstein SMMS implies an upper bound for the diameter and a positive lower bound for the volume on the underlying manifold, both of which are uniform in $m\geq 1+\delta$.

5) The tensor $A$ appearing in the fifth assumption is the weighted Weyl tensor, which we define in Section~\ref{sec:generalized_curvature}.  As we shall see in Corollary~\ref{cor:rm_equiv_a}, for compact quasi-Einstein SMMS with \charconstant\ one and $m\geq 1+\delta>1$, a uniform $L^{n/2}$-bound for $A$ is equivalent to a uniform $L^{n/2}$-bound for the Riemann curvature tensor $\Rm$.  We have chosen to use the weighted Weyl curvature due to its importance in establishing the $\varepsilon$-regularity lemma and in understanding the orbifold singularities.  With an eye towards generalizing our precompactness theorem, it is important to note that $A$ does not carry enough information to prove Theorem~\ref{thm:precpt} in the cases $m=0,1$ (see Section~\ref{sec:generalized_curvature}).

6) The computations in Section~\ref{sec:generalized_curvature} will show that the uniform bound on the $L^{n/2}$ norm of $A$ is equivalent to a uniform bound on the $L^{n/2}$ bound of the Riemann curvature tensor in the quasi-Einstein scale.  Furthermore, in Section~\ref{sec:geometry}, we will see that there is a uniform bound $\lv\nabla f_i\rv\leq C(n,\delta,C_1,C_2)$ for the quasi-Einstein scales.  Since the remaining assumptions of Theorem~\ref{thm:precpt} are conformally invariant, this shows that the statement and conclusion of Theorem~\ref{thm:precpt} also hold in the quasi-Einstein scale.

7) If $n=3$, one expects to be able to classify quasi-Einstein SMMS, analogous to the situation for Einstein metrics (they are spaceforms) and shrinking gradient Ricci solitons (they are rigid in the sense of~\cite{Petersen_Wylie1}; see~\cite{Zhang2008b}).  In particular, one would have a considerably stronger conclusion than that of Theorem~\ref{thm:precpt} for sequences of three-dimensional quasi-Einstein SMMS.  For a partial result in this direction, see~\cite{CMMR2010}.

8) As we will see in Section~\ref{sec:proof}, the orbifold singularities of $M$ arise as regions where the sectional curvatures of $M_i$ blow up.  By our assumptions, we will see that this can only happen at isolated points, and moreover, that a sequence of suitable rescalings of the metrics $g_i$ in annuli around a point of curvature concentration will converge smoothly to metric cones $C(S^{n-1}/\Gamma)$ for $\Gamma\subset\mathrm{O}(n)$ a finite group acting freely on $S^{n-1}$.  The group $\Gamma$ is then easily seen to correspond to the orbifold group of the corresponding singular point, and moreover, it will be seen that both the number $\#\mS$ of singular points and the size $|\Gamma|$ of $\Gamma$ are bounded above in terms of $n,\delta,C_1,C_2$.  Additionally, since there are no nontrivial (i.e.\ $|\Gamma|>1)$ groups $\Gamma\subset\mathrm{O}(n)$ for $n$ odd such that $S^{n-1}/\Gamma$ is orientable, we see that if $n$ is odd and the $M_i$ are assumed to be orientable, the limiting space $M$ is in fact smooth.

9) It is important to note that our definition of convergence does \emph{not} rule out the possible formation of bubbles; indeed, there are many examples of convergent sequences of Einstein metrics which do bubble off topology (see, for example, \cite{AndersonCheeger1991}).  In particular, as the above discussion suggests, our result only requires that we know that the ``necks'' connecting bubbling regions to $M_i$ (or other bubbling regions) look like annuli in the cones $C(S^{n-1}/\Gamma)$ centered at the vertex.

10) Bando~\cite{Bando1990} and, from a different perspective, Anderson and Cheeger~\cite{AndersonCheeger1991} have carried out a more detailed study of the aforementioned bubbling behavior in the context of convergent sequences of Einstein metrics.  In particular, they were able to show that bubbles carry a definite amount of energy $\int\lv\Rm\rv^{n/2}\geq\theta>0$ for $\theta$ dependent only on $n,\delta,C_1,C_2$, and Anderson and Cheeger further were able to establish finer control on the geometry of the necks to conclude that there are only finitely many diffeomorphism types of Einstein metrics which meet the assumptions of Theorem~\ref{thm:precpt}.  The author expects that it is possible to extend their result to our setting, but as the proof does not seem to require additional insights into the geometric consequences of positive $m$-energy, has opted not to pursue this question here.

\vspace{10pt}

The proof of Theorem~\ref{thm:precpt} will proceed in a series of steps, which are standard in proving results of this type.  The key step is applying the Cheeger-Gromov precompactness theorem~\cite{Cheeger1970,Gromov1981} to a collection of open balls which almost cover the manifolds $M_i$, from which standard arguments allow us to prove our theorem.  In order to apply this theorem, we need to establish uniform bounds on the volume and sectional curvature of each ball.  To find the convergent subsequence, we also need global control of the volume and the diameter of the manifolds.  Then, in order to prove that the convergence is smooth, we need to be able to control higher derivatives of the curvatures and potentials.  Because these arguments are rather standard, we shall only sketch the proof, highlighting the changes that must be made in our setting (cf.\ \cite{Anderson1989,Cao_Sesum,Weber2008,Zhang2006,Zhang2008}).

To avoid tedious indexing, we will repeatedly use the symbols $C, C_i$, $i\geq 3$ to denote bounds on various quantities related to the SMMS of Theorem~\ref{thm:precpt} which are independent of $m$ and $\lambda$, but which are allowed to change from line to line.

\subsection{The Weighted Weyl Curvature}
\label{sec:generalized_curvature}

Though the weighted Weyl curvature $A$ appearing in Theorem~\ref{thm:precpt} is in some sense equivalent to the Riemann curvature tensor, it has the benefit that it is divergence-free with respect to a weighted measure and that it satisfies a nice elliptic equation (see Theorem~\ref{thm:div_free_riem} and Theorem~\ref{thm:lapl_riem_conformal}, respectively).  These properties are necessary in understanding the orbifold singularities and in establishing the $\varepsilon$-regularity lemma, respectively.  For an explanation of the notation used in this section, as well as for a derivation of the results presented, we refer the reader to the appendix.

\begin{defn}
Let $(M^n,g,1^m\dvol)$ be a SMMS with \charconstant\ $\mu$.  The \emph{weighted Weyl curvature $A$} is defined by
\begin{equation}
\label{eqn:A}
A = \Rm - P\wedge g ,
\end{equation}
where
\begin{equation}
\label{eqn:P}
P = \frac{1}{m+n-2}\left(\Ric - \frac{R+m\mu}{2(m+n-1)}g\right)
\end{equation}
and $\wedge$ denotes the Kulkarni-Nomizu product of symmetric $(0,2)$-tensors.
\end{defn}

\begin{remark}
The weighted Weyl curvature first appeared in the case $m<2-n$ in work of He, Petersen and Wylie~\cite{HePetersenWylie2010}, who used the fact that it is divergence-free in a suitable weighted sense when $(M^n,g,1^m\dvol)$ admits a quasi-Einstein scale to study the local properties of locally conformally flat metrics.
\end{remark}

In the case $m=0$, it is clear that $A$ is the Weyl tensor and $P$ is the Schouten tensor; it is for this reason that we have opted to call $A$ the weighted Weyl tensor.  Note also that when $m=\infty$, the weighted Weyl tensor is the Riemannian curvature tensor, and thus the assumptions of Theorem~\ref{thm:precpt} coincide with the usual assumption in theorems of this type (cf.\ \cite{Weber2008,Zhang2006,Zhang2009}).

\begin{remark}
If we insist when rescaling SMMS that the measure is always the usual Riemannian measure, we must also rescale the \charconstant\ to ensure that geometric objects associated to a SMMS with \charconstant\ scale properly.  More precisely, given a SMMS $(M^n,g,1^m\dvol_g)$ with \charconstant\ $\mu$, the constant $c>0$ determines the rescaled SMMS $(M^n,c^{2}g,1^m\dvol_{c^{2}g})$ with \charconstant\ $c^{-2}\mu$.  In this way, the weighted Weyl curvature, regarded as a section of $S^2\Lambda^2T^\ast M$, rescales as
\[ A\left(c^{2}g,1^m\dvol_{c^{2}g},c^{-2}\mu\right) = c^2 A\left(g,1^m\dvol_g,\mu\right) . \]
\end{remark}

It is clear that $A$ is an algebraic curvature tensor, and so we can use the orthogonal decomposition
\[ \lp C \rp = \lp W \rp \oplus \lp \Ric_0 \rp \oplus \lp \id \rp \]
of algebraic curvature tensors into operators of Weyl-type, traceless-Ricci-type, and multiples of the identity (cf.\ \cite{Besse,BohmWilking2006}) to rewrite $A$ in the following useful form:

\begin{prop}
Let $(M^n,g,1^m\dvol)$ be a SMMS with \charconstant\ $\mu$.  Then the weighted Weyl curvature can be written
\begin{equation}
\label{eqn:a_decomposed}
\begin{split}
A & = W + \frac{m}{(m+n-2)(n-2)}\Ric_0 \wedge g \\
& \quad + \frac{m}{2(m+n-1)(m+n-2)}\left(\frac{m-1}{n(n-1)}R + \mu\right)g\wedge g .
\end{split}
\end{equation}
where $W$ is the Weyl curvature and $\Ric_0=\Ric-\frac{1}{n}Rg$ is the traceless part of the Ricci curvature.
\end{prop}

\begin{proof}

This follows immediately from the definition of $A$ and the usual Ricci decomposition
\begin{equation}
\label{eqn:weyl_decomposition}
\Rm = W + \frac{1}{n-2}\Ric_0\wedge g + \frac{1}{2n(n-1)}Rg\wedge g . \qedhere
\end{equation}
\end{proof}

\begin{remark}

We note that, when $m\not=0,1$, the weighted Weyl curvature vanishes if and only if $(M^n,g)$ is a spaceform with constant sectional curvature $-\frac{\mu}{m-1}$.  In particular, for the SMMS considered in Theorem~\ref{thm:precpt}, the only examples for which $A\equiv0$ are the positive elliptic Gaussians~\eqref{eqn:hyperbolic_peg}.
\end{remark}

As a consequence of this decomposition, it is straightforward to compare the norms of $A$ and $\Rm$.

\begin{prop}
\label{prop:rm_norm_a}
Let $(M^n,g,1^m\dvol)$ be a SMMS with \charconstant\ one, and suppose additionally that $m\geq 1+\delta$ for some $\delta>0$.  Then there is a constant $1<C<\infty$ depending only on $n$ and $\delta$ such that
\[ \left|A\right|^2 \leq \left|\Rm+\frac{1}{2(m-1)}g\wedge g\right|^2 \leq C\left|A\right|^2 . \]
\end{prop}

\begin{remark}
$-\frac{1}{2(m-1)}g\wedge g$ is the Riemann curvature tensor of a spaceform $\overline{M}$ with constant sectional curvature $-(m-1)^{-1}$, and so this says that $A$ is small if and only if the sectional curvature of $M$ is close to that of $\overline{M}$.  Since $m-1\geq\delta>0$, the curvature of $\overline{M}$ is itself uniformly bounded, and it is in this sense that we will use Proposition~\ref{prop:rm_norm_a}.
\end{remark}

\begin{proof}

Since the Ricci decomposition~\eqref{eqn:weyl_decomposition} is orthogonal, we have that
\begin{align}
\notag \left|\Rm+kg\wedge g\right|^2 & = \left|W\right|^2 + \frac{1}{n-2}\left|\Ric_0\right|^2 + \frac{1}{2n(n-1)}(R+2n(n-1)k)^2 \\
\label{eqn:norm} \left|A\right|^2 & = \left|W\right|^2 + \frac{1}{n-2}\left(\frac{m}{m+n-2}\right)^2\left|\Ric_0\right|^2 \\
\notag & \quad + \frac{n(n-1)}{2}\left(\frac{m}{(m+n-1)(m+n-2)}\right)^2\left(\frac{m-1}{n(n-1)}R+1\right)^2
\end{align}
for any $k\in\bR$.  In particular, using $k=\frac{1}{2(m-1)}$ and the assumption $m>1$, it holds that
\[ \left|A\right|^2\leq \left|\Rm+\frac{1}{2(m-1)}g\wedge g\right|^2 \leq \left(\frac{(m+n-1)(m+n-2)}{m(m-1)}\right)^2\left|A\right|^2 . \]
The result follows by observing that the coefficient on the right is at least one and is bounded above in terms of $n$ and $\delta$.
\end{proof}

This allows us to verify our claim that the uniform $L^{n/2}$-bound on $A$ is equivalent to a uniform $L^{n/2}$-bound on $\Rm$ in Theorem~\ref{thm:precpt}.  A key ingredient in verifying this claim, which is also necessary to establish diameter bounds for the SMMS in Theorem~\ref{thm:precpt}, is the following estimate on the volume of such a SMMS in terms of $L^{n/2}$-bounds on $A$ and $\Rm$ (cf.\ \cite{Zhang2009}).

\begin{lem}
\label{lem:vol_upper_bound}
Let $(M^n,g,1^m\dvol)$ be a SMMS with \charconstant\ one which admits a quasi-Einstein scale $u\in C^\infty(M)$.  Then there are positive constants $K_1,K_2$ depending only on $n$ such that
\begin{align*}
\Vol(M) & \leq K_1 \left\|A\right\|_{n/2}^{n/2} \\
\Vol(M) & \leq K_2 \left\|\Rm\right\|_{n/2}^{n/2} .
\end{align*}
\end{lem}

\begin{proof}

By Corollary~\ref{cor:duality}, $(M^n,g,u^{2-m-n}\dvol)$ is a quasi-Einstein SMMS with quasi-Einstein constant one.  Taking the trace of~\eqref{eqn:qe_ric} and integrating  with respect to $1^m\dvol$, we see that
\[ \int_M \left(R + \frac{1}{m+n-2}|\nabla f|^2\right)\dvol = n\Vol(M) . \]
On the other hand, integrating~\eqref{eqn:bianchi} and using the fact that the quasi-Einstein constant $\lambda$ is positive yields
\[ \frac{1}{m+n-2}\int_M |\nabla f|^2 < \Vol(M) . \]
Combining these two integrals, we see that
\begin{align*}
\int_M R\,\dvol & > (n-1)\Vol(M) \\
\int_M \left(\frac{m-1}{n(n-1)}R + 1\right)\dvol & > \frac{m+n-1}{n}\Vol(M) .
\end{align*}
It then follows from H\"older's inequality that
\begin{align*}
\Vol(M)^{\frac{2}{n}} & \leq \frac{1}{n-1}\left(\int_M R^{\frac{n}{2}}\dvol\right)^{\frac{2}{n}} \\
\Vol(M)^{\frac{2}{n}} & \leq \frac{n}{m+n-1} \left(\int_M \left(\frac{m-1}{n(n-1)}R+1\right)^{\frac{n}{2}}\dvol\right)^{\frac{2}{n}} .
\end{align*}
Using~\eqref{eqn:norm}, we then see that
\begin{align*}
\Vol(M)^{\frac{2}{n}} & \leq \sqrt{\frac{2n}{n-1}} \left(\int_M \left|\Rm\right|^{\frac{n}{2}}\dvol\right)^{\frac{2}{n}} \\
\Vol(M)^{\frac{2}{n}} & \leq \frac{m+n-2}{m}\sqrt{\frac{2n}{n-1}} \left(\int_M \left| A\right|^{\frac{n}{2}}\dvol\right)^{\frac{2}{n}} .
\end{align*}
Since $\frac{m+n-2}{m}<n-1$ when $m>1$, this yields the explicit constants $K_1,K_2$.
\end{proof}

\begin{cor}
\label{cor:rm_equiv_a}
Given $C_1,C_2>0$ and $n\geq 3$, there are positive constants $K_1(n,C_1)$, $K_2(n,C_2)$ such that whenever $(M^n,g,1^m\dvol)$ is a SMMS with \charconstant\ one and $m>1$ which admits a quasi-Einstein scale, then
\begin{align*}
\left\|A\right\|_{n/2} \leq C_1 \Longrightarrow \left\|\Rm\right\|_{n/2} \leq K_1 \\
\left\|\Rm\right\|_{n/2} \leq C_2 \Longrightarrow \left\|A\right\|_{n/2} \leq K_2 .
\end{align*}
\end{cor}

\begin{proof}

This follows directly from Proposition~\ref{prop:rm_norm_a} and Lemma~\ref{lem:vol_upper_bound}.
\end{proof}

\begin{remark}
Arguing in a similar way and using the transformation formula for the Schouten tensor for a conformal change of metric (see~\cite{Besse}), one can show that if $(M^n,g,v^m\dvol_g)$ is a quasi-Einstein SMMS with $m\geq1+\delta$, \charconstant\ one, and $\lV\Rm\rV_{n/2}\leq C$, then the SMMS
\[ \left( M^n, v^{-2}g, 1^m\dvol_{v^{-2}g}\right) \]
with \charconstant\ one is such that $\lV A\rV_{n/2}\leq C^\prime$, where $C^\prime$ depends only on $n,C,\delta$.
\end{remark}

As we mentioned, the reason for introducing the weighted Weyl curvature is the following nice expressions for its divergence and its Laplacian, which are derived in the appendix.

\begin{thm}
\label{thm:div_free_riem}
Let $(M^n,g,1^m\dvol)$ be a SMMS with \charconstant\ $\mu$ which admits a quasi-Einstein scale $u=e^{\frac{f}{m+n-2}}$.  Then the weighted Weyl curvature satisfies
\[ \delta A = \frac{m+n-3}{m+n-2}\imath_{\nabla f} A . \]
\end{thm}

\begin{thm}
\label{thm:lapl_riem_conformal}
Let $(M^n,g,1^m\dvol)$ be a SMMS with \charconstant\ $\mu$ which admits a quasi-Einstein scale $u=e^{\frac{f}{m+n-2}}$.  Then the weighted Weyl curvature satisfies
\begin{equation}
\label{eqn:lapl_riem_conformal}
\Delta_f A = 2\mu A - A^2 - A^\# - \frac{1}{m}\tr A\wedge\tr A + \frac{2}{(m+n-2)^2}\lp A,df\otimes df\rp \wedge g ,
\end{equation}
where $\Delta_f A = \Delta f - \imath_{\nabla f}A$.
\end{thm}

For our purposes, all we really need from Theorem~\ref{thm:lapl_riem_conformal} is that
\begin{equation}
\label{eqn:lapl_riem_conformal_ast}
\Delta_f A = 2\mu A + A\ast A + A\ast df^2,
\end{equation}
where $A\ast B$ denotes a tensor constructed by taking linear combinations of tensors formed from $A\otimes B$ by switching indices and contracting with the metric; in particular, using the exact form in which the coefficients depend on $m\geq1+\delta>1$, we know that $\lv A\ast A\rv\leq C\lv A\rv^2$ and $\lv A\ast df^2\rv\leq C\lv A\rv\,\lv df\rv^2$ for some constant $C$ depending only on $n$ and $\delta$.
\subsection{Controlling Diameter and Volume}
\label{sec:geometry}

In order to apply the Cheeger-Gromov theorem to the quasi-Einstein SMMS of Theorem~\ref{thm:precpt}, we need control of their diameters as well as both local and global volume estimates.  We have already established uniform upper and lower bounds on the volume of the SMMS through Lemma~\ref{lem:vol_upper_bound} and Proposition~\ref{prop:global_volume_bound_qe}, respectively, while Theorem~\ref{thm:kappa_qe} yields uniform lower bounds $\Vol(B(x,r))\geq\kappa r^{-n}$ on the volumes of small balls.  To get the remaining estimates, we need the following growth estimate for the scalar curvature, the quasi-Einstein scale, and its derivative:

\begin{prop}
\label{prop:bounds_no_m}
Let $(M^n,g,1^m\dvol)$ be a compact SMMS with \charconstant\ $\mu>0$ which admits a quasi-Einstein scale $u=e^{\frac{f}{m+n-2}}$ with quasi-Einstein constant $\lambda$.  Then for all $x\in M$, the estimate
\begin{equation}
\label{eqn:quad_growth}
\sup \left\{ R(x), |\nabla f|^2(x), f(x) \right\} \leq C_3 r(x)^2 + C_4
\end{equation}
holds, where $r(x)$ is the distance $d(x,p)$ for $p$ a minimizer of $f$ and $C_3,C_4$ depend on $m,n,\mu,\lambda,f(p)$.

Moreover, if $(M^n,g,1^m\dvol)$ and $u$ are normalized as in Theorem~\ref{thm:precpt}, then $C_1,C_2$ depend only on $n,\delta,C_1,C_2$.
\end{prop}

\begin{proof}

In terms of $f$, \eqref{eqn:dil_estimate} is
\[ \lv\nabla f\rv^2 < \frac{(m+n-2)^2}{m-1}\mu - \frac{(m+n-2)^2}{m+n-1}\lambda e^{-\frac{2}{m+n-2}f} . \]
Using the elementary estimate $e^{-ax}\geq 1-ax$, this yields
\begin{equation}
\label{eqn:linearized_dil}
\lv\nabla f\rv^2 \leq 2\frac{m+n-2}{m+n-1}\lambda f + (m+n-2)^2\left(\frac{\mu}{m-1}-\frac{\lambda}{m+n-1}\right) =: k_1 f + k_2 .
\end{equation}
It thus follows that $f\leq k_3 r^2 + k_4$ for constants $k_3>0$ and $k_4$ depending on $n,\delta,\lambda,m(\mu-\lambda),\inf f$.  Plugging back in to~\eqref{eqn:linearized_dil} yields a similar bound for $\lv\nabla f\rv^2$, while~\eqref{eqn:bianchi} and~\eqref{eqn:qe_uv} imply that
\begin{align*}
R & = -\frac{m+n-1}{m+n-2}\lv\nabla f\rv^2 - (m+n-2)\lambda e^{-\frac{2}{m+n-2}f} + (m+2n-2)\mu \\
& \leq 2\lambda f + (m+2n-2)\mu-(m+n-2)\lambda,
\end{align*}
yielding~\eqref{eqn:quad_growth}.

The final claim will follow once we can establish uniform upper bounds on $\lambda$, $\inf f$, and $m(1-\lambda)$.  As a first observation, Corollary~\ref{prop:compute_energy_quasi_einstein}, Proposition~\ref{prop:global_volume_bound_qe} and Lemma~\ref{lem:vol_upper_bound} imply that there are constants $0<a_1,a_2$ depending only on $n,\delta,C_1,C_2$ such that $a_1\leq\lambda^{\frac{m+n}{2}}\leq\Vol(M)\leq a_2$.  In particular, $\lambda\leq\max\{a_2^{2/n},1\}$ and
\[ m(1-\lambda) \leq m\left(1-a_1^{\frac{2}{m+n}}\right) . \]
Since $m\mapsto m(1-a_1^{\frac{2}{m+n}})$ is increasing in $m$ when $a_1\leq 1$, this implies that $m(1-\lambda)\leq\max\{-2\log a_1,0\}$.  Finally, \eqref{eqn:qe_lambda_uv} implies that
\[ \inf f\leq\frac{m+n-2}{2}\log \lambda \leq \log a_2, \]
yielding the necessary bound on $\inf f$.
\end{proof}

Combining Proposition~\ref{prop:compute_energy_quasi_einstein}, Proposition~\ref{prop:global_volume_bound_qe}, Theorem~\ref{thm:vol_growth_log}, Lemma~\ref{lem:vol_upper_bound}, and Proposition~\ref{prop:bounds_no_m}, we thus have the following result which includes most of the non-curvature estimates we need to prove Theorem~\ref{thm:precpt}.

\begin{prop}
\label{prop:diameter_bound}
Let $(M^n,g,1^m\dvol)$ be a SMMS satisfying the hypotheses of Theorem~\ref{thm:precpt}.  Then there is a constant $C(n,\delta,C_1,C_2)>0$ such that
\[ \diam(M) \leq C, \quad \lv\nabla f\rv^2\leq C, \quad \Vol(M) \geq C^{-1} . \]
\end{prop}

The final estimate we need is an upper bound on the growth of the volume of balls, which is a consequence of a Bishop-Gromov-type volume comparison theorem.

\begin{prop}
\label{prop:vol_comp}
Let $(M^n,g,1^m\dvol)$ be a SMMS satisfying the hypotheses of Theorem~\ref{thm:precpt}.  Then there is a constant $C_3(n,\delta,C_1,C_2)>0$ such that for all $0<r_1<r_2<\diam(M)$,
\[ \frac{\Vol(B(x,r_2))}{\Vol(B(x,r_1))}\leq C_3\left(\frac{r_2}{r_1}\right)^n . \]
\end{prop}

\begin{proof}

By Corollary~\ref{cor:dil}, it holds that $\frac{1}{m+n-2}df^2\leq\frac{m+n-2}{m-1}$g, which can be uniformly bounded in terms of $n$ and $\delta$.  Using Corollary~\ref{cor:duality}, it follows that there is a $c_1(n,\delta)\geq 0$ such that $\Ric + \nabla^2f\geq -c_1g$.  By Proposition~\ref{prop:diameter_bound}, there are uniform bounds on $|\nabla f|^2$, and so we may apply the volume comparison theorem of Wei and Wylie~\cite[Theorem 1.1(b)]{Wei_Wylie} to yield the result.
\end{proof}
\subsection{An $\varepsilon$-Regularity Lemma}
\label{sec:epsilon_regularity}

In order to apply the Cheeger-Gromov theorem, we also need to establish uniform estimates for the sectional curvature.  These estimates come in the form of a so-called $\varepsilon$-regularity lemma, which provides uniform estimates on the sectional curvatures of balls for which the $L^{n/2}$ norm of the Riemann curvature tensor $\Rm$ is small.  One can establish such estimates provided $\Rm$ satisfies a suitably nice elliptic equation and the manifold is $\kappa$-noncollapsed (cf.\ \cite{TianViaclovsky2008,Weber2008}).  It is for this reason that we have introduced the weighted Weyl curvature, as Theorem~\ref{thm:lapl_riem_conformal} provides such an elliptic equation.

To make the above precise, we shall follow the arguments of Weber~\cite{Weber2008} which yield the $\varepsilon$-regularity lemma for compact shrinking gradient Ricci solitons.  As the modifications necessary to adapt his results to our setting are relatively minor, we shall only point out the necessary modifications and refer to~\cite{Weber2008} for the remaining details.

\begin{lem}
\label{lem:weber}
Let $(M^n,g,1^m\dvol)$ be a SMMS with \charconstant\ one and $m\geq 1+\delta>1$.  Then there are constants $\varepsilon(n,\delta),C(n,\delta)>0$ such that for all balls $B=B(p,r)$ on which $\sup_B \left|A\right|\leq\varepsilon$,
\[ C_s(B) \leq C \left(\frac{\Vol(B(p,r))}{r^n}\right)^{-2/n} , \]
where $C_s(B)$ is the Sobolev constant on $B$.
\end{lem}

\begin{proof}

By Proposition~\ref{prop:rm_norm_a}, $\lv A\rv\leq\varepsilon$ if and only if $\Rm$ is close to the curvature of the spaceform with constant sectional curvature $-(m-1)^{-1}\geq-\delta^{-1}$.  In particular, we may choose $\varepsilon=\varepsilon(n,\delta)$ small enough so that $\exp_p$ is noncritical on $B(0,1)\subset T_pM$ for all $p\in M$ and moreover, such that any ball of radius one which is diffeomorphic to a Euclidean ball and has $\lv A\rv\leq\varepsilon$ must have Sobolev constant bounded above by $2C_E$, where $C_E$ is the Euclidean Sobolev constant.  On such a ball, the Bishop volume comparison estimate
\[ \Vol(B(p,1)) \leq 2v(n,-\delta^{-1},1) = c(n,\delta) \]
holds, where $v(n,H,r)$ is the volume of the $n$-dimensional spaceform of constant sectional curvature $H$.  With these ingredients, one can follow Weber's argument~\cite[Lemma~4.1]{Weber2008}.
\end{proof}

\begin{prop}[{\cite[Proposition~4.2]{Weber2008}}]
\label{prop:weber}
Let $(M^n,g,1^m\dvol)$ be a SMMS with \charconstant\ one and $m\geq 1+\delta>1$.  Suppose additionally that $\varepsilon$-regularity holds for $A$; that is, suppose that for all $B(r)=B(x,r)$,
\[ C_s^{\frac{n}{2}} \int_{B(r)} \left|A\right|^{\frac{n}{2}}\dvol \leq \varepsilon \Longrightarrow \sup_{B(r/2)} \left|A\right| \leq C(1+r^{-2})\left(C_s^{\frac{n}{2}}\int_{B(r)}\left|A\right|^{\frac{n}{2}}\dvol\right)^{\frac{2}{n}} \]
for constants $\varepsilon,C$ depending only on $n$.  Then there are (possibly different from above) constants $\varepsilon(n,\delta),C(n,\delta)>0$ such that, if
\[ H = \sup_{B(q,s)\subset B(p,r)} \frac{s^n}{\Vol(B(q,s))}\int_{B(q,s)} \left|A\right|^{\frac{n}{2}}\dvol \]
satisfies $H\leq\varepsilon$, then
\[ \sup_{B(p,r/2)} \left|A\right| \leq C(1+r^{-2})H^{\frac{2}{n}} . \]
\end{prop}

\begin{proof}

Weber's proof~\cite[Proposition~4.2]{Weber2008} carries through verbatim, except that one must use Lemma~\ref{lem:weber} in place of Weber's version for $\Rm$.
\end{proof}

Together, these two results yield the following $\varepsilon$-regularity theorem, which does not depend on Sobolev constant $C_s$.

\begin{thm}[$\varepsilon$-regularity Lemma]
\label{thm:eps_regularity}
Let $(M^n,g,1^m\dvol)$ be a SMMS satisfying the hypotheses of Theorem~\ref{thm:precpt}.  Then there are constants $C,\varepsilon>0$ depending only on $n,\delta,C_1$ such that if
\[ \int_{B(x,r)} \lv A\rv^{n/2}1^m\dvol \leq \varepsilon, \]
then
\[ \sup_{B(x,r/2)} \lv A\rv \leq C(1+r^{-2})\left(\int_{B(x,r)} \lv A\rv^{n/2}\dvol\right)^{2/n} . \]
\end{thm}

\begin{proof}

As discussed above, we know that $\Delta_f A = c_4A + A\ast A + A\ast df^2$, where the constants depend only on $n$ and $\delta$.  Hence
\begin{align*}
\lp A, \Delta A \rp & = \lp A, \nabla_{\nabla f} A\rp + \lp A, \Delta_f A \rp \\
& \geq -\frac{1}{4}\lv\nabla A\rv^2 - (c_4+\lv\nabla f\rv^2)|A|^2 - c_5\lv A\rv^3,
\end{align*}
where the constants $c_4,c_5$ depend only on $n$ and $\delta$.  On the other hand, using Kato's inequality $\left|\nabla\lv A\rv\right|^2\leq\lv\nabla A\rv^2$, we have
\[ \lv A\rv\Delta\lv A\rv = \frac{1}{2}\Delta\lv A\rv^2 - \left|\nabla\lv A\rv\right|^2 \geq \lp A,\Delta A\rp . \]
By Proposition~\ref{prop:diameter_bound}, there is a constant $c_3(n,C_1,C_2)>0$ such that $\lv df\rv^2\leq c_3$.  Combining these estimates, it then follows that
\[ \lv A\rv\Delta\lv A\rv + \frac{1}{4}\lv \nabla A\rv^2 + c_6\lv A\rv^2 + c_5\lv A\rv^3 \geq 0 . \]
By Moser iteration (cf.\ \cite{Anderson1989,Weber2008}), it then follows that there are constants $\varepsilon(n,c_5,c_6)>0$, $c_7(n,c_5,c_6)>0$ such that whenever $B(x,r)$ is a geodesic ball such that
\[ C_s^{\frac{n}{2}}\int_{B(x,r)} \lv A\rv^{\frac{n}{2}} \dvol \leq \varepsilon, \]
then
\[ \sup_{B(x,r/2)} \lv A\rv \leq c_7(1+r^{-2})\left(C_s^{\frac{n}{2}} \int_{B(x,r/2)}\lv A\rv\dvol\right)^{\frac{2}{n}}, \]
where $C_s$ is the Sobolev constant of $(M,g)$.  By Theorem~\ref{thm:kappa_qe}, $(M,g)$ is $\kappa$-noncollapsed, and thus Proposition~\ref{prop:weber} yields the conclusion.
\end{proof}
\subsection{Shi-type Estimates for Higher Derivatives}
\label{sec:shi_estimates}

In order to establish smooth convergence, we also need estimates on the derivatives of the curvature (cf.\ \cite{Shi1989}).  The necessary estimates are all local estimates, and so it is beneficial to find well-adapted cutoff functions.  To do so, the following Laplacian comparison lemma is useful:

\begin{lem}
\label{lem:lapl_comp}
Let $(M^n,g,1^m\dvol)$ be a SMMS satisfying the hypotheses of Theorem~\ref{thm:precpt}.  Given a point $p\in M$, let $r$ denote the distance from $p$, $r(x)=d(p,x)$.  Then there are constants $C_3,C_4>0$ depending only on $n,\delta,C_1,C_2$ such that
\[ \Delta_f r \leq \frac{C_3}{r}(1+C_4 r) \]
holds in the barrier sense.
\end{lem}

\begin{proof}

As in the proof of Proposition~\ref{prop:vol_comp}, there is a constant $c_1(n,\delta,C_1,C_2)>0$ such that $\Ric+\nabla^2f\geq -c_1 g$.  The result then follows from the Laplacian comparison theorem of Wei and Wylie~\cite[Theorem 2.1]{Wei_Wylie} together with the elementary estimate $\coth(cr)\leq 1+\frac{1}{cr}$.
\end{proof}

As a consequence, we can construct the desired cutoff functions:

\begin{lem}
\label{lem:cutoff_fn}
Let $(M^n,g,1^m\dvol)$ be a complete SMMS with \charconstant\ one and $m\geq 1+\delta>1$ which admits a quasi-Einstein scale $u=e^{\frac{f}{m+n-2}}$.  Let $U\subset M$ be an open set containing the closure of the ball $B(x,a)$ for some $a\leq 1$.  Then there is a function $\eta\colon U\to[0,1]$ such that
\begin{enumerate}
\item $C_3 a^2\leq \eta\leq 1$ on $B(x,a/2)$
\item $\eta\equiv 0$ on $U-B(x,a)$
\item $\lv\nabla\eta\rv^2 \leq C_4$
\item $\Delta_f\eta \geq -C_5$ ,
\end{enumerate}
where (4) holds in the barrier sense and $C_3,C_4,C_5$ are positive constants depending only on $n,\delta,C_1,C_2$.
\end{lem}

\begin{proof}

One can take the function $\eta(y) = \eta(d(x,y))$, where $\eta\colon\bR\to[0,1]$ is defined by $\eta(r)=a^2-r^2$.  Since $a\leq 1$, the desired properties follow immediately from Lemma~\ref{lem:lapl_comp}.
\end{proof}

Using these cutoff functions, we establish the desired estimates on higher derivatives of the curvature.

\begin{thm}
\label{thm:shi}
Let $(M^n,g,1^m\dvol)$ be a SMMS as in Theorem~\ref{thm:precpt}.  Let $U\subset M$ be an open set containing the closure of the ball $B(x,a)$ for $a\leq 1$, and suppose
\[ \sup_U \left|\Rm\right| \leq C_3 . \]
Then for all $k\geq 0$,
\[ \sup_{B(x,a/2)} \left(\left|\nabla^k\Rm\right| + \left|\nabla^{k+2}f\right|\right) \leq \frac{C(n,\delta,C_1,C_2,C_3)}{a^k} . \]
\end{thm}

\begin{proof}

Using the cutoff functions of Lemma~\ref{lem:cutoff_fn}, we will in fact prove the slightly stronger estimate
\[ \sup_{U} \eta^k\left(\lv\nabla^k\Rm\rv^2+\lv\nabla^{k+2}f\rv^2\right) \leq C(n,\delta,C_1,C_2,C_3,k) , \]
which we shall establish by a relatively straightforward induction argument.

We first note that it suffices to show $\eta^k\left|\nabla^kA\right|^2\leq C(n,\delta,C_1,C_2,C_3,k)$.  To see this, observe that by the definition of $A$,
\[ \left|\nabla^k\Rm\right|^2 \leq C(n,k,\delta) \left|\nabla^k A\right|^2 . \]
By Corollary~\ref{cor:duality},
\[ \frac{1}{3}\left|\nabla^2 f\right|^2 \leq n + \left|\Ric\right|^2 + \left(\frac{1}{m+n-2}\right)^2 \left|df\right|^4 , \]
whence $\left|\nabla^2 f\right|\leq C$.  By taking derivatives, we then see that $\left|\nabla^{k+2}f\right|$ is controlled by $\left|\nabla^k\Rm\right|$ and $\left|\nabla^{k+1}f\right|$.  Thus, by induction, we see that $\left|\nabla^{k+2}f\right|$ is controlled by $\left|\nabla^i\Rm\right|$ for $0\leq i\leq k$.

Now, by Theorem~\ref{thm:lapl_riem_conformal} and the Ricci identity, we find that
\begin{align*}
\Delta_f \left|\nabla^kA\right|^2 & \geq c(n,\delta)\left|\nabla^{k+1}A\right|^2 - c_1\left|\nabla^kA\right|^2 - \sum_{i=0}^k c_3(i) \left|\nabla^kA\right| \left|\nabla^{k-i}A\right| \left|\nabla^{i}A\right| \\
& \quad - \sum_{i+j=0}^{k-1} c_4(i,j) \left|\nabla^kA\right| \left|\nabla^{i}A\right| \left|\nabla^{j+1}f\right| \left|\nabla^{k-i-j-1}A\right| \\
& \quad - \sum_{i=0}^k c_5(i) \left|\nabla^kA\right| \left|\nabla^{i+1}f\right| \left|\nabla^{k-i}A\right| ,
\end{align*}
where the constants $c_i$ depend only on $n$ and $\delta$.  Let $\eta$ be a cutoff function as in Lemma~\ref{lem:cutoff_fn}.  Using the Schwarz inequality
\[ 2\eta^{k+1}\left|\nabla^{k}A\right| \left|\nabla^{i}A\right| \left|\nabla^{k-i}A\right| \leq \eta^{k}\left|\nabla^kA\right|^2 + \left(\eta\left|\eta^{i/2}\nabla^iA\right| \left|\eta^{(k-i)/2}\nabla^{k-i}A\right|\right)^2, \]
one can easily show that
\begin{align*}
\Delta_f \left|\eta^{\frac{k+1}{2}}\nabla^kA\right|^2 & = \eta^{k+1}\Delta_f\left|\nabla^kA\right|^2 + \left|\nabla^kA\right|^2\Delta_f\eta^{k+1} \\
& \quad + 4(k+1)\lp\left|\nabla^kA\right|\eta^{\frac{k-1}{2}}\nabla\eta,\eta^{\frac{k+1}{2}}\nabla\left|\nabla^kA\right|\rp \\
& \geq 2\eta^{k+1}\left|\nabla^{k+1}A\right|^2 - c_6\eta^k\left|\nabla^kA\right|^2 - c_7,
\end{align*}
where $c_6,c_7>0$ depend only $n,k,\delta,C_1,C_2,C_3$.  Thus we can find positive constants $\kappa_i$ depending only on $n,k,\delta,C_1,C_2,C_3$ such that
\[ E = \sum_{i=0}^k \kappa_i \eta^i \left|\nabla^iA\right|^2 \]
satisfies
\[ \Delta_f E \geq E - c_8  \]
for some constant $c_8$ depending only on $n,k,\delta,C_1,C_2,C_3$.  Hence the maximum principle yields $E\leq c_8$, and in particular,
\[ \eta^k\left|\nabla^kA\right|^2 \leq c_8, \]
as desired.
\end{proof}

\begin{remark}
By tracking the dependence of the constants $C(n,\delta,C_1,C_2,C_3,k)$ on $C_3$, we see that if we instead assume that
\[ \sup_U \eta^l|\Rm|^2\leq C_3, \]
then we can conclude that there are constants $C(n,\delta,C_1,C_2,C_3,k)$ such that
\[ \sup_U \eta^{k+l}\left(|\nabla^k\Rm|^2 + |\nabla^{k+2}f|^2\right) \leq C . \]
We shall in fact use Theorem~\ref{thm:shi} in this form with $l=2$, as we will start with the initial estimate $\lv\Rm\rv\leq C_3 r^{-2}$ from Proposition~\ref{prop:rm_norm_a} and Theorem~\ref{thm:eps_regularity}.
\end{remark}
\subsection{Proof of the Compactness Theorem}
\label{sec:proof}

We now turn to the proof of Theorem~\ref{thm:precpt}.  As these arguments are standard, we only sketch the argument; additional details can be found in~\cite{Anderson1989,Weber2008,Zhang2006,Zhang2009}.

\begin{proof}[Proof of Theorem~\ref{thm:precpt}]

First observe that, using the estimates of Proposition~\ref{prop:diameter_bound} and Proposition~\ref{prop:vol_comp}, Gromov's theorem~\cite{Gromov1981} implies that $(M_i,g_i)$ converges in the Gromov-Hausdorff topology to $(M,g)$, which is \emph{a priori} only a length space.

Next, let $\varepsilon$ be as in the $\varepsilon$-regularity theorem and fix $r>0$ small.  Let $x_i^k$ be a maximal $r/4$-separated set of points in $M_i$, so that $M_i=\bigcup B(x_i^k, r/2)$ and $B(x_i^k,r/8)$ are disjoint.  Define the sets of good and bad balls by
\begin{align*}
G_{i,r} & = \left\{ B(x_i^k,r/2) \colon \int_{B(x_i^k,r)} \left|A_i\right|^{n/2}\dvol_i \leq \varepsilon \right\} \\
B_{i,r} & = \left\{ B(x_i^k,r/2) \colon \int_{B(x_i^k,r)} \left|A_i\right|^{n/2}\dvol_i \geq \varepsilon \right\} ,
\end{align*}
respectively.  Since $\varepsilon$ depends only on $n,\delta,C_1,C_2$, and since Proposition~\ref{prop:vol_comp} yields a uniform upper bound on the number of disjoint balls of radius $r/8$ inside a ball of radius $2r$ in $M$, we have a uniform upper bound on the number of bad balls $\# B_{i,r}$ depending only on $n,\delta,C_1,C_2$.  Moreover, by passing to a subsequence, we may assume that $Q=\# B_{i,r}$ is constant.  By Theorem~\ref{thm:eps_regularity}, we have uniform bounds on $\left|A\right|$ inside the good balls $B(x_i^k,r/2)\subset G_{i,r}$, and hence, by Proposition~\ref{prop:rm_norm_a}, uniform bounds on $\left|\Rm\right|$.  Theorem~\ref{thm:shi} then yields uniform estimates on $|\nabla^j\Rm(g_i)|$ and $|\nabla^j f_i|$ for all $i,j$.  Letting $G_r^i=\bigcup B(x_i^k,r/2)$, where the union is only over those balls in $G_r^i$, we then see that there is a smooth quasi-Einstein SMMS $G_r$ such that $G_r^i\to G_r$ in the Cheeger-Gromov sense, after possibly passing to a subsequence.

Next, choose a sequence $r_{j+1}<\frac{r_j}{2}\to 0$ and find a subsequence $G_{r_j}$ of quasi-Einstein SMMS which are smooth limits as above.  Set $G=\bigcup G_{r_j}$ with the induced metric from each $G_{r_j}$.  By the above discussion, this is a smooth quasi-Einstein SMMS.  Let $\overline{G}$ be the metric completion of $G$.  As in~\cite{Anderson1989}, one can check that $\overline{G}=G\cup\mS$, for $\mS=\{p_1,\dotsc,p_k\}$ a set of isolated points (in $\overline{G}$) with $k\leq Q$, and moreover, that $\overline{G}=M$.

To show that $M$ is an orbifold, it thus remains to establish that singular points $p\in\mS$ are orbifold points.  Since we have established uniform $C^1$ bounds on $f$ on $M\setminus\mS$, this essentially reduces to the problem of ensuring that the metric $g$ can be extended to a $C^2$ orbifold metric, as the smoothness of $g$ and $f$ then follow from elliptic regularity (cf.\ Section~\ref{sec:shi_estimates}).  Establishing this regularity proceeds as in~\cite{Anderson1989}, with some modifications necessary for our setting.  More precisely, we may assume without loss of generality that $\mS=\{p\}$, and establishing that $p$ is an orbifold point proceeds in three steps:

First, we show that $p$ contains a punctured neighborhood for which each connected component is diffeomorphic to a metric cone $C(S^{n-1}/\Gamma)$.  To see this, consider annular regions $A(s/l,ls)=\{x\in M\colon s/l\leq d(x,p)\leq ls\}$ for $s$ small.  Theorem~\ref{thm:kappa_qe} and Proposition~\ref{prop:vol_comp} yield uniform constants $\kappa,c$ such that $\kappa r^{n}\leq\Vol(B(y,r))\leq cr^{n}$ for all $B(y,r)\subset A(s/l,ls)$, while Theorem~\ref{thm:eps_regularity} and Theorem~\ref{thm:shi} imply that $|\nabla^k\Rm|^2(x)\leq c_k d(x,p)^{-k-2}$ for uniform constants $c_k$.  Since these estimates are all scale invariant, we can take a diagonal limit of the rescaled annuli $(A(s/l,ls),s^{-2}g)$ as $s\to 0$ and $l\to \infty$ which converges to the metric cone $C(S^n/\Gamma)$.  Moreover, as in~\cite{Anderson1989}, the $\kappa$-noncollapsing estimate $\Vol(B(y,r))\geq\kappa r^{n}$ gives a uniform bound $|\Gamma|\leq N$ by comparing volume growth in the annuli and in $C(S^n/\Gamma)$.

Second, we can show that $p$ is a $C^0$ orbifold point.  This follows by showing that for $n$ large and $s$ small (in sense depending only on $n,\delta,C_1,C_2$), at most one component of $A(s/l,ls)$ intersects $\partial B(p,s)$, which in turn follows directly from the arguments in~\cite{Anderson1989} by using Proposition~\ref{prop:vol_comp} in place of the Bishop-Gromov volume comparison (cf.\ \cite{Zhang2009}).  Thus $p$ has a punctured neighborhood $U$ which looks like $C(S^n/\Gamma)$.  Moreover, as in~\cite{Anderson1989}, the fact that $A(s/2,2s)$ converges to $C(S^n/\Gamma)$ smoothly on compact sets allows us to conclude that the metric $g$ extends continuously to the universal cover of $U$, and hence $p$ is a $C^0$ orbifold point.

Finally, we can show that $p$ is in fact a $C^\infty$ orbifold point by showing that the metric constructed above is in fact smooth.  This follows from elliptic regularity theory once we know that $\int_U |A|^q<\infty$ for some $q>n/2$.  On the other hand, we do know that $\int_U |A|^{n/2}<\infty$, and so we can adapt the well-known techniques for removing singularities of Yang-Mills fields~\cite{Sibner1985,Uhlenbeck1982a} to achieve the desired $L^q$ bound.  More precisely, the results of Section~\ref{sec:epsilon_regularity} imply that the Sobolev constant is uniformly bounded below on $M\setminus\mS$.  Thus, if $n\geq 5$, we can use argument of Sibner~\cite{Sibner1985} to show that $|A|\in L^q(U)$ for all $q$, while if $n=4$, we can modify Uhlenbeck's argument~\cite{Uhlenbeck1982a} by using Theorem~\ref{thm:div_free_riem} and the uniform bounds on $|\nabla f|^2$ from Proposition~\ref{prop:diameter_bound} in place of the Yang-Mills condition $\delta\Rm=0$ (cf.\ \cite{Weber2008,Zhang2006,Zhang2009}).
\end{proof}                 % Precompactness

\appendix
\section{The Weighted Laplacian of $A$}
\label{sec:computation}

The purpose of this appendix is to prove Theorem~\ref{thm:div_free_riem} and Theorem~\ref{thm:lapl_riem_conformal}.  It is possible to carry out these computations by directly taking derivatives of the quasi-Einstein equation together with the definition of curvature.  However, this is a rather cumbersome process, largely due to the large number of terms which appear due to the presence of the quadratic term $df^2$.  To avoid this, we will find it useful to first discuss the Weitzenb\"ock formula relating the rough Laplacian to the Hodge Laplacian on $S^2\Lambda^2T^\ast M$.

\subsection{The Weitzenb\"ock Formula}
\label{sec:computation/weitzenbock}

First, we introduce some notation.  Given a Riemannian manifold $(M^n,g)$, the Levi-Civita connection determines canonical operators
\begin{align*}
d^1&\colon\Lambda^kT^\ast M\otimes\Lambda^lT^\ast M \to \Lambda^{k+1}T^\ast M\otimes\Lambda^lT^\ast M \\
\delta^1&\colon\Lambda^kT^\ast M\otimes\Lambda^lT^\ast M\to\Lambda^{k-1}T^\ast M\otimes\Lambda^lT^\ast M \\
d^2&\colon\Lambda^kT^\ast M\otimes\Lambda^lT^\ast M\to\Lambda^kT^\ast M\otimes\Lambda^{l+1}T^\ast M\\
\delta^2&\colon\Lambda^kT^\ast M\otimes\Lambda^lT^\ast M\to\Lambda^kT^\ast M\otimes\Lambda^{l-1}T^\ast M
\end{align*}
for any $0\leq k,l\leq n$, which are the twisted exterior derivative and its divergence (with respect to $\dvol_g$) taken on the first and second factors, respectively.  When computing, we will always fix a point $p\in M$ and let $x,y,z\in T_pM$ denote vector fields which are evaluated in the first factor and $u,v$ be vector fields which are evaluated on the second factor.  Our conventions in defining $d^i$ and $\delta^i$ are such that for $T\in T^\ast M\otimes T^\ast M$ and $A\in\Lambda^2T^\ast M\otimes\Lambda^2T^\ast M$,
\begin{align*}
(d^1T)(x,y,u) & = \nabla_x T(y,u) - \nabla_x T(y,u) \\
(\delta^1T)(u) & = \sum_{i=1}^n \nabla_{e_i} T(e_i,u) \\
(d^1A)(x,y,z,u,v) & = \nabla_x A(y,z,u,v) + \nabla_y A(z,x,u,v) + \nabla_z A(x,y,u,v) \\
(\delta^1A)(y,u,v) & = \sum_{i=1}^n \nabla_{e_i} A(e_i,y,u,v),
\end{align*}
where $\{e_i\}$ is an orthonormal basis of $T_pM$ and all vector fields have been extended to a neighborhood of $p$ by parallel translation.  Note here that we are thinking of $T$ and $A$ as sections of the respective vector bundles, an abuse of notation we shall make throughout this appendix.  Also, we will adopt the convention that the Riemann curvature tensor $\Rm\in S^2\Lambda^2T_p^\ast M$ is defined by
\[ \Rm(x,y,u,v) := \lp -\nabla_x\nabla_yu+\nabla_y\nabla_xu, v \rp =: \lp R(x,y)u,v\rp ; \]
recall in this formula that we have specified that $\nabla x=0=\nabla y$ at $p$.  With this convention, the Ricci curvature $\Ric\in S^2T^\ast M$ is defined by
\[ \Ric(x,u) = \sum_{i=1}^n \Rm(e_i,x,e_i,u) . \]

Our interest is in relating the rough Laplacian $\Delta$ to the Hodge Laplacian $\Delta_H$ on (sections of) $S^2\Lambda^2T^\ast M$, which are defined by
\begin{align*}
(\Delta A)(x,y,u,v) & = \nabla_{e_i}\nabla_{e_i}A(x,y,u,v) \\
\Delta_H A & = \frac{1}{2}\left(\delta^1d^1+d^1\delta^1+\delta^2d^2+\delta^2d^2\right)A ,
\end{align*}
where we have adopted Einstein summation notation.  In fact, we will be interested in the weighted analogues $\Delta_\phi$ and $\Delta_{\phi,H}$ defined in the obvious way given a measure $e^{-\phi}\dvol_g$ on $(M^n,g)$.

In order to carry out our computation, we also need to consider a number of different algebraic operators.  First, as a general operator on the bundles $\Lambda^kT^\ast M\otimes\Lambda^lT^\ast M$ for various values of $k$ and $l$, we have the \emph{wedge product}
\[ \wedge \colon \left(\Lambda^kT^\ast M\otimes\Lambda^lT^\ast M\right)\times\left(\Lambda^rT^\ast M\otimes\Lambda^sT^\ast M\right) \to \Lambda^{k+r}T^\ast M\otimes\Lambda^{l+s}T^\ast M \]
is defined in the obvious way on each factor.  For example, the Kulkarni-Nomizu product is just the wedge product of two sections $h,k\in T^\ast M\otimes T^\ast M$:
\[ (h\wedge k)(x,y,u,v) = h(x,u)k(y,v) + h(y,v)k(x,u) - h(x,v)k(y,u) - h(y,u)k(x,v) . \]

\begin{remark}
This convention is such that $g\wedge g=2\id$, where we use the metric $g$ to regard $\id\in\End(\Lambda^2T^\ast M)$ as a section of $S^2\Lambda^2T^\ast M$ via the natural inner product for which the two-form $e_i\wedge e_j$ has length one for $i\not=j$.
\end{remark}

There are two actions involving $S^2T^\ast M$ and $S^2\Lambda^2T^\ast M$ which we shall need.  First is the contraction
\[ \lp\cdot,\cdot\rp\colon S^2\Lambda^2T^\ast M \times S^2T^\ast M \to S^2T^\ast M \]
defined by
\[ \lp A,T\rp(x,u) = A\left(e_i,x,T(e_i),u\right) . \]
In particular, $\lp A,g\rp$ is just the trace of $A$; e.g.\ $\Ric=\lp\Rm,g\rp$.  Second is the action
\[ \hash \colon S^2T^\ast M\times S^2\Lambda^2T^\ast M\to S^2\Lambda^2T^\ast M \]
defined by
\begin{align*}
(T\hash A)(x,y,u,v) & = -A\left(T(x),y,u,v\right) - A\left(x,T(y),u,v\right) \\
& \quad - A\left(x,y,T(u),v\right) - A\left(x,y,u,T(v)\right) .
\end{align*}
In other words, $T\hash$ is the natural action of $T\in T^\ast M\otimes T^\ast M\cong\End(TM)$ on $TM$ extended as a derivation to the tensor algebra of $TM$.

On $S^2\Lambda^2T^\ast M$ there are two natural notions of a product.  First, by using the metric to identify $S^2\Lambda^2T^\ast M$ with the space of symmetric endomorphisms of $\Lambda^2T^\ast M$, one defines the product
\[ \circ \colon S^2\Lambda^2T^\ast M \times S^2\Lambda^2T^\ast M \to \Lambda^2T^\ast M\otimes \Lambda^2T^\ast M \]
by
\begin{equation}
\label{eqn:composition_product}
(B\circ A)(x,y,u,v) = A(x,y,e_i,e_j) B(e_i,e_j,u,v) .
\end{equation}

\begin{remark}
As defined above, $(B\circ A)(x,y,u,v) = 2\lp B\left(A(x\wedge y)\right),u\wedge v\rp$.  The convention~\eqref{eqn:composition_product} is more common in the literature, which is our reason for adopting it.
\end{remark}

Of course, \eqref{eqn:composition_product} does not in general produce sections of $S^2\Lambda^2T^\ast M$, and for this reason we will more typically consider the symmetric product
\[ A\cdot B = \frac{1}{2}\left( A\circ B + B\circ A \right) . \]

The second product is the natural Lie algebra product arising from the (fiber-wise) identification $\Lambda^2T^\ast M\cong\kso(n)$.  Explicitly, one defines the Lie bracket $[\cdot,\cdot]$ on $\Lambda^2T^\ast M$ by
\[ [\alpha,\beta](x,y) = \lp\imath_x\alpha,\imath_y\beta\rp - \lp\imath_y\alpha,\imath_x\beta\rp \]
to realize $\Lambda^2T^\ast M$ as a Lie algebra.  One then defines the product
\[ \sq \colon S^2\Lambda^2T^\ast M\times S^2\Lambda^2T^\ast M \to S^2\Lambda^2T^\ast M \]
by
\[ (\alpha\otimes\alpha)\sq(\beta\otimes\beta) = [\alpha,\beta]\sq[\alpha,\beta] \]
and extending by linearity (cf.\ \cite{BohmWilking2006,Hamilton1986}).  For general $A,B\in S^2\Lambda^2T^\ast M$, this yields the formula
\begin{align*}
(A\sq B)(x,y,u,v) & = A(x,e_i,u,e_j) B(y,e_i,v,e_j) + A(y,e_i,v,e_j)B(x,e_i,u,e_j) \\
& \quad - A(x,e_i,v,e_j)B(y,e_i,u,e_j) - A(y,e_i,u,e_j)B(x,e_i,v,e_j) .
\end{align*}
Following~\cite{Hamilton1986}, we will also abbreviate $A^\#=A\sq A$.

As a consequence, we are now in a position to state and prove the Weitzenb\"ock formula relating $\Delta_\phi$ and $\Delta_{\phi,H}$.

\begin{thm}
\label{thm:weitzenbock}
Let $(M^n,g,e^{-\phi}\dvol_g,m)$ be a SMMS and let $A\in S^2\Lambda^2T^\ast M$.  Then it holds that
\begin{equation}
\label{eqn:weitzenbock}
\Delta_\phi A = \Delta_{\phi,H}A - \frac{1}{2}\Ric_\phi^\infty\hash A - \Rm\cdot A - \Rm\sq A .
\end{equation}
\end{thm}

\begin{proof}

First, suppose that $\phi$ is a constant.  Using the Ricci identity, it is straightforward to verify that
\begin{align*}
(\Delta A)(x,y,u,v) & = \nabla_{e_i}\nabla_{e_i}A(x,y,u,v) \\
& = \nabla_{e_i}d^1A(e_i,x,y,u,v) + \nabla_{e_i}\nabla_x A(e_i,y,u,v) - \nabla_{e_i}\nabla_y A(e_i,x,u,v) \\
& = \left((\delta^1d^1+d^1\delta^1)A\right)(x,y,u,v) + A\left(\Ric(x),y,u,v\right) + A\left(x,\Ric(y),u,v\right) \\
& \quad + A\left(e_i,R(e_i,x)y,u,v\right) + A\left(e_i,R(y,e_i)x,u,v\right) + A\left(e_i,y,R(e_i,x)u,v\right) \\
& \quad + A\left(e_i,y,u,R(e_i,x)v\right) - A\left(e_i,x,R(e_i,y)u,v\right) - A\left(e_i,x,u,R(e_i,y)v\right) \\
& = \left((\delta^1d^1+d^1\delta^1)A\right)(x,y,u,v) + A\left(\Ric(x),y,u,v\right) + A\left(x,\Ric(y),u,v\right) \\
& \quad - (A\circ\Rm)(x,y,u,v) - (A\sq\Rm)(x,y,u,v) .
\end{align*}
Symmetrizing then yields
\[ \Delta A = \Delta_H A - \frac{1}{2}\Ric\hash A - \Rm\cdot A - \Rm\sq A . \]

Next, let $\phi$ be arbitrary and observe that
\begin{align*}
(\nabla_{\nabla\phi} A)(x,y,u,v) & = d^1A(\nabla\phi,x,y,u,v) + \nabla_xA(\nabla\phi,y,u,v) - \nabla_yA(\nabla\phi,x,u,v) \\
& = d^1A(\nabla\phi,x,y,u,v) + d^1\left(A(\nabla\phi)\right)(x,y,u,v) \\
& \quad - A\left(\nabla_x\nabla\phi,y,u,v\right) - A\left(x,\nabla_y\nabla\phi,u,v\right) ,
\end{align*}
where we have written $A(\nabla\phi):=A(\nabla\phi,\cdot,\cdot,\cdot)$, a convention we shall employ for the remainder of this appendix.  Symmetrizing this and subtracting it from the previous display then yields the desired result.
\end{proof}

The last ingredient we need is the following simple computational lemma, whose proof we shall omit.

\begin{lem}
\label{lem:trivial_facts_ap}
Given $\alpha\in\Lambda^kT^\ast M\otimes T^\ast M$ and $f\in C^\infty(M)$, it holds that
\begin{align*}
d^1(\alpha\wedge g) & = d^1\alpha\wedge g \\
\delta^1(\alpha\wedge g) & = \delta^1\alpha\wedge g + (-1)^{k+1}d^2\alpha \\
\imath_{\nabla f}(\alpha\wedge g) & = \imath_{\nabla f}\alpha\wedge g + (-1)^k \alpha\wedge df,
\end{align*}
where the contraction $\imath_{\nabla f}\colon\Lambda^{k+1}T^\ast M\otimes\Lambda^2T^\ast M\to\Lambda^kT^\ast M\otimes\Lambda^2T^\ast M$ is taken in the first factor and we regard $df\in\Lambda^0T^\ast M\otimes\Lambda^1T^\ast M$ in the final equation.
\end{lem}
\subsection{Proofs of Theorem~\ref{thm:div_free_riem} and Theorem~\ref{thm:lapl_riem_conformal}}
\label{sec:computation/weyl}

Let us now turn to the verification of Theorem~\ref{thm:div_free_riem} and Theorem~\ref{thm:lapl_riem_conformal}.  To that end, recall that we are studying a SMMS $(M^n,g,1^m\dvol)$ with \charconstant\ $\mu$ which admits a quasi-Einstein scale $u=e^{\frac{f}{m+n-2}}$.  This assumption on $u$ in particular forces
\begin{equation}
\label{eqn:qe_appendix}
\Ric + \nabla^2 f + \frac{1}{m+n-2}df\otimes df = \mu g,
\end{equation}
which is the main fact we will need.  Our objective is then to find formulae for the divergence and the Laplacian of the weighted Weyl curvature $A$, which we recall is defined by
\begin{align*}
A & = \Rm - P \wedge g \\
P & = \frac{1}{m+n-2}\left(\Ric - \frac{R+m\mu}{2(m+n-1)}g\right) .
\end{align*}

To start, using nothing but the definitions of $A$ and $P$ and the Bianchi identity, the following facts are easily verified.

\begin{lem}
\label{lem:basic_facts_ap}
With $A$ and $P$ as above,
\begin{align}
\label{eqn:dP_genl} dP & = \frac{1}{m+n-2}\left(d\Ric - \frac{1}{2(m+n-1)}dR\wedge g\right) \\
\label{eqn:divP_genl} \delta P & = \frac{1}{2(m+n-1)}dR \\
\label{eqn:AdP_genl} \delta^2 A & = (m+n-3)d^1P .
\end{align}
\end{lem}

\begin{proof}

\eqref{eqn:dP_genl} follows immediately from the definition of $P$, \eqref{eqn:divP_genl} from the Bianchi identity $\delta\Ric=\frac{1}{2}dR$, and~\eqref{eqn:AdP_genl} follows from~\eqref{eqn:divP_genl} and the Bianchi identity $\delta^2\Rm=d^1\Ric$.
\end{proof}

Theorem~\ref{thm:div_free_riem} is then an immediate consequence of the following result.

\begin{prop}
\label{prop:div_free_app}
Let $(M^n,g,1^m\dvol_g)$ be a SMMS with \charconstant\ $\mu$ and quasi-Einstein scale $u$ as above.  Then
\begin{equation}
\label{eqn:Af_qe}
A(\nabla f) = (m+n-2)dP .
\end{equation}
In particular, Theorem~\ref{thm:div_free_riem} holds.
\end{prop}

\begin{remark}
Here and in the following, we shall simply write $d$ in place of $d^1$ or $d^2$ whenever the meaning is clear from context.
\end{remark}

\begin{proof}

From Lemma~\ref{lem:trivial_facts_ap} and the definition of $A$, it holds in general that
\[ A(\nabla f) = \Rm(\nabla f) - P(\nabla f)\wedge g - df\wedge P . \]
Now, on the one hand, applying $d^1$ to~\eqref{eqn:qe_appendix} implies that
\begin{align*}
\Rm(\nabla f) & = (m+n-2)dP + df\wedge P \\
& \quad + \frac{(m+n-2)dR + R\,df - (m+2n-2)\mu\,df}{2(m+n-1)(m+n-2)}\wedge g .
\end{align*}
On the other hand, contracting $\nabla f$ into~\eqref{eqn:qe_appendix} implies that
\begin{align*}
-(m+n-2)P(\nabla f) & = \frac{1}{2}d\lv\nabla f\rv^2 + \frac{1}{m+n-2}\lv\nabla f\rv^2\,df \\
& \quad + \frac{1}{2(m+n-1)}\left(R-(m+2n-2)\mu\right)df .
\end{align*}
Hence, it follows that
\begin{align*}
A(\nabla f) & = (m+n-2)dP \\
& \quad + \frac{u^{-2}}{2(m+n-1)}d\left(u^2(R+\frac{m+n-1}{m+n-2}\lv\nabla f\rv^2 - (m+2n-2)\mu)\right)\wedge g \\
& = (m+n-2)dP,
\end{align*}
where the final equality follows from~\eqref{eqn:bianchi}.

Finally, combining~\eqref{eqn:AdP_genl} and~\eqref{eqn:Af_qe} yields the proof of Theorem~\ref{thm:div_free_riem}.
\end{proof}

As an immediate consequence of the Weitzenb\"ock formula~\eqref{eqn:weitzenbock}, we have the following result.

\begin{cor}
\label{cor:almost_lapl_riem}
Let $(M^n,g,1^m\dvol_g)$ be a SMMS with \charconstant\ $\mu$ and quasi-Einstein scale $u$ as above.  Then
\begin{align*}
\Delta_f A & = 2\mu A - \Rm\cdot A - \Rm\sq A \\
& \quad - \frac{1}{m+n-2}\left(\lp A,\nabla^2f-\frac{1}{m+n-2}df\otimes df\rp\wedge g\right) .
\end{align*}
\end{cor}

\begin{proof}

By Theorem~\ref{thm:div_free_riem}, we have that
\begin{equation}
\label{eqn:useful_identity_app}
\delta_f A = -\frac{1}{m+n-2}A(\nabla f) ,
\end{equation}
while the Bianchi identity $d\Rm=0$ and Proposition~\ref{prop:div_free_app} together imply that
\[ d^2A = \delta_f^1 A \wedge g . \]
Thus we may write $dA=\delta_f A\wedge g$, whence follows
\[ \delta_f^2d^2A = \delta_f^2\delta_f^1A\wedge g - d^1\delta_f^1A + \delta_f^1 A\wedge df . \]
Applying~\eqref{eqn:useful_identity_app}, we see that
\[ \delta_f^2\delta_f^1A = -\frac{1}{m+n-2}\lp A,\nabla^2f-\frac{1}{m+n-2}df\otimes df\rp . \]
Thus it follows that
\[ (\delta_f^2d^2+d^1\delta_f^1)A = -\frac{1}{m+n-2}\left(\lp A, \nabla^2f-\frac{1}{m+n-2}df\otimes df\rp\wedge g + A(\nabla f)\wedge df\right) . \]

On the other hand, it is straightforward to check that
\[ \left(A(\nabla f)\wedge df\right)(x,y,u,v) = -A\left(x,y,df^2(u),v\right) - A\left(x,y,u,df^2(v)\right) . \]
After symmetrizing, the result then follows from~\eqref{eqn:weitzenbock} and~\eqref{eqn:qe_appendix}.
\end{proof}

To prove Theorem~\ref{thm:lapl_riem_conformal}, it thus remains to verify the following algebra lemma.

\begin{lem}
With $A=\Rm-P\wedge g$, it holds that
\begin{equation}
\label{eqn:algebra_lemma}
A^2 + A^\# = \Rm\cdot A + \Rm\sq A - \tr A\wedge P - \lp A,P\rp\wedge g .
\end{equation}
\end{lem}

\begin{proof}

Clearly we have that
\[ A^2 + A^\# = \Rm\cdot A + \Rm\sq A - A\cdot(P\wedge g) - A\sq(P\wedge g) . \]
On the other hand, it is straightforward to check that
\begin{align*}
\left(A\circ(P\wedge g)\right)(x,y,u,v) & = 2A\left(P(x),y,u,v\right) + 2A\left(x,P(y),u,v\right) \\
A\sq(P\wedge g) & = \tr A\wedge P + \lp A,P\rp\wedge g + P\hash A ,
\end{align*}
from which the result immediately follows.
\end{proof}

This allows us to complete the proof of Theorem~\ref{thm:lapl_riem_conformal}.

\begin{proof}[Proof of Theorem~\ref{thm:lapl_riem_conformal}]

From Corollary~\ref{cor:almost_lapl_riem} and~\eqref{eqn:algebra_lemma} it follows that
\begin{equation}
\label{eqn:app_proof1}
\begin{split}
\Delta_f A & = 2\mu A - A^2 - A^\# - \tr A\wedge P \\
& \quad - \lp A, P + \frac{1}{m+n-2}\nabla^2 f - \frac{1}{(m+n-2)^2}df^2\rp\wedge g .
\end{split}
\end{equation}
From~\eqref{eqn:qe_appendix} we see that
\begin{align}
\notag P & + \frac{1}{m+n-2}\nabla^2 f - \frac{1}{(m+n-2)^2}df^2 \\
\label{eqn:app_proof2}& = -\frac{1}{m+n-2}\left(\frac{2}{m+n-2}df^2 + \frac{R-(m+2n-2)\mu}{2(m+n-1)}g\right) .
\end{align}
On the other hand, it is straightforward to check from the definitions of $A$ and $P$ that
\begin{equation}
\label{eqn:app_proof3}
\tr A = m\left(P - \frac{R-(m+2n-2)\mu}{2(m+n-1)}g\right) .
\end{equation}
Combining~\eqref{eqn:app_proof1}, \eqref{eqn:app_proof2}, and~\eqref{eqn:app_proof3} together with the fact $\lp A,g\rp=\tr A$ then yields~\eqref{eqn:lapl_riem_conformal}, as desired.
\end{proof}            % The messy \Delta A computation

\bibliographystyle{abbrv}
\bibliography{../bib}
\end{document}